\documentclass[amssymb,12pt ]{article}
\usepackage{geometry} 
\usepackage{graphicx} 
\usepackage{amsfonts}
\usepackage{amssymb}
\usepackage{amsmath}

\usepackage{longtable}

\def\Tri{\text{\rm Tri }}
\def\Area{\text{\rm Area }}
\def\vol{\text{\rm vol }}

\def\al {\aligned}
\def\eal {\endaligned}

\def\vv {\Vert}
\def\la {\langle}
\def\ra {\rangle}
\newtheorem{thm}{Theorem}[section]
\newtheorem{cor}[thm]{Corollary}
\newtheorem{lem}[thm]{Lemma}
\newtheorem{pro}[thm]{Proposition}

\newtheorem{defi}[thm]{Definition}

\newtheorem{rem}[thm]{Remark}
\newtheorem{exa}[thm]{Example}

\newenvironment{proof}{{\bf Proof }}{\hfill\framebox[2mm]{}}

\begin{document}

\title{The analytic continuation of hyperbolic space}
\author{Yunhi Cho and Hyuk Kim\footnote{The second author was supported by grant no.(R01-1999-000-00002-0)
  from the Basic Research Program of the Korea Science $\&$ Engineering Foundation}}
\date{}
\maketitle

\begin{abstract}
We define and study an extended hyperbolic space which contains
the hyperbolic space and de Sitter space as subspaces and which is
obtained as an analytic continuation of the hyperbolic space. The
construction of the extended space gives rise to a complex valued
geometry consistent with both the hyperbolic and de Sitter space.
Such a construction shed a light and inspires a new insight for
the study of the hyperbolic geometry and Lorentzian geometry. We
discuss the advantages of this new geometric model as well as some
of its applications.

\end{abstract}

\renewcommand \figurename{Fig}
\makeatletter
\renewcommand\fnum@figure[1]{\textbf{\figurename} }
\makeatother
\section{Introduction}

The hyperbolic space is an independent geometric entity with an
infinite diameter and infinite volume which is already complete in
its own right. But if we look at the hyperbolic space as a unit disk
in the Kleinian model, then using the same metric formula we have a
Lorentzian space with constant curvature outside the unit disk.
Furthermore we can even draw a geometric figure lying across the
ideal boundary. We naturally expect on this space the generalization
of the basic geometric notions such as angle, length, volumes, ...,
etc, and the similar relation between them to those on the
hyperbolic space. But we immediately have difficulties in defining
and deriving those due to the Lorentzian nature of the metric and
multi-valuedness of the analytic functions representing various
geometric formulas.

In this paper we show there is a natural way of extending the
geometry of hyperbolic space to Lorentzian part and set a
foundation of a geometry  which connects and unifies these two
different geometries by an analytic continuation method on the
Kleinian model. We call such a unified space an extended
hyperbolic space since we start from the hyperbolic space and then
continue analytically to the Lorentzian part. The purpose of
studying an extended hyperbolic space is not only to give a proper
geometry on the Lorentzian space as a continuation of hyperbolic
geometry but also to give a new insight to the hyperbolic geometry
itself by studying a geometric object lying across the ideal
boundary - the proper study of such object would be impossible
otherwise.

Some of the basic notions such as angle, length and geodesic
triangle on the extended hyperbolic plane has been considered and
studied so far through cross ratio (see \cite{Sc}) and a rather ad
hoc combinatorial method (see \cite{D}). In this paper we study the
extended hyperbolic space in a more systematic way and discuss some
of its applications. In Section 2, we set up the geometry as an
analytic continuation going over the singularity of the hyperbolic
metric at the ideal boundary, we view the hyperbolic metric as a
limit of complex perturbation called an $\epsilon$-metric which is a
complex regular metric. And we define and study distance, angle,
length and k-dimensional volumes on the extended space in Section 4.
Then such geometric quantities are given rather naturally with
complex numbers - they of course coincide with the usual real values
for the quantities inside the hyperbolic space. If we consider a
nice region lying across the ideal boundary, then the volume of the
region is of a finite complex value while the volume of the
hyperbolic part and the Lorentzian part of the region are both
infinite. In fact the measure defined on the extended hyperbolic
space which extends the usual hyperbolic volume (and Lorentz volume)
is a finitely (but not countably) additive complex measure and we
study some of its delicate and strange properties in Section 3.

We would like to describe some of the advantages and expectations
of using the extended hyperbolic model. When we study hyperbolic
geometry we naturally want to extend the geometric objects and
notions over and beyond the ideal boundary. But when we try to
compute geometric quantities, we come across with a confusion of
choosing an appropriate value among infinite possibilities of
multi-values even for a distance or an angle. The extended model
provide us a fulfilling consistency without such confusion once we
choose an analytic continuation which follows naturally after
fixing an $\epsilon$-approximation of the metric. For instance the
understanding of 1-dimensional extended model quickly leads us to
be able to define an angle as a complex numbers on the general
semi-Riemannian manifolds. Generalizing an earlier work of
Kellerhals\cite{K}, Ushijima studied the volume of a hyperbolic
polyhedron obtained by truncating with the dual plane of an ultra
ideal vertex of a tetrahedron which is lying across the ideal
boundary. He showed in \cite{U} that the volume of the truncated
tetrahedron is the real part of the value obtained after formal
application of the known volume formula of a hyperbolic
tetrahedron to this tetrahedron. But for the imaginary part it is
multi-valued and there remains a problem of choice and
interpretation of its geometric meaning. The extended model
determines the unique value of the imaginary part and explains its
geometric meaning as the volume of the truncated portion or as the
area of the truncated face. See Example 5.11.

As we work more with the extended model, we found that the model
is a natural and fundamental geometric setting as it works
beautifully in every aspect. The n-dimensional extended hyperbolic
space $\Bbb S^n_H$ is simply the standard sphere $\Bbb S^n$
topologically but with the geometry coming from the unit sphere of
the Minkowski space. Its geometry has many resemblance with the
spherical geometry and we can obtain various results through this
conceptual but concrete analogy. For instance we can derive the
Gauss-Bonnet formula on the extended space (this was first shown
in \cite{I}) combinatorially using Euler method as we did for the
2-sphere $\Bbb S^2$ without computing integrals and then of course
can be generalized to the higher dimensions. (See Proposition
\ref{4.2}.) Furthermore we can extend the hyperbolic trigonometry
to the extended space, which of course implies that we have the
same trigonometry on the de Sitter space (the Lorentzian part) and
even for an object sitting across the ideal boundary. In principle
we can apply the same method used for the standard sphere to
obtain the corresponding results on the extended hyperbolic space,
which in turn give rise to the results for the hyperbolic space as
well as de Sitter space. It would be interesting to observe that
the volume of $\Bbb S^n_H$ differs from that of $\Bbb S^n$ by
$i^n$. (See Theorem 2.3.)

As another illustration, let us consider Santal\'o's formula
\cite{S} giving the relation between the volume of a simplex and
its dual in $\Bbb S^3$ or $\Bbb H^3$, Milnor's relation \cite{M}
between a convex polyhedron $P$ and its dual $P^*$ in $\Bbb S^3$,
and Su\'arez-Peir\'o's result \cite{E} for a simplex and its dual
in $\Bbb H^n$. Then we see that there is a slight discrepancy
between $\Bbb S^3$ and $\Bbb H^3$:
 $$\al &\text{vol }(P)+\text{vol }(P^*)+\frac
12\sum a_i(\pi-A_i)=\pi^2\quad\quad ~(\Bbb S^3)\\&\text{vol
}(P)+\text{vol }(P^*)-\frac 12\sum a_i(\pi-A_i)=0\quad\quad
~~(\Bbb H^3)\eal
$$
where $a_i$ and $A_i$ are edge lenghts and dihedral angles of $P$.
Here if we use the extended model $\Bbb S^3_H$ instead of $\Bbb
H^3$, the proof is essentially identical with the case for $\Bbb
S^3$ and the formula can be written for both cases in a unified
way as follows.
$$\text{vol }(P)-\text{vol
}(\rm II)+\text{vol }(\rm III)-\text{vol }(\rm IV)+\text{vol
}(P^*)=0$$ For this expression, the volume of simplices of type
{\rm III} lying across the boundary play the crucial role which is
also interpreted as mean curvature integral in $(\Bbb S^3)$ and
$(\Bbb H^3)$. And Su\'arez-Peir\'o's result is then also expected
to be generalized for both cases in a single identity as
$$\text{vol }(P)-\text{vol }({\rm II})+\text{vol }({\rm
III})-\cdots+(-1)^{n+1}\text{vol }(P^*)=0$$
 Conceptually things
are getting tremendously easier giving us inspirations and this is
one of the real merits using the extended model.

Lastly here are some speculations and the directions for possible
further developments and applications. If we use the extended
model, all the geometric quantities become complex valued and we
can ask whether this is related to complex invariants of
hyperbolic manifolds such as volume and Chern-Simon invariant pair
in dimension 3. Using the extended model we can derive all the
trigonometry and precise elementary geometric formulas for the
Lorentzian spherical space (or de Sitter space) and this will be
useful to study the discrete group actions on the Lorentzian
space. Since the hyperbolic space and the Lorentzian space are
dual each other, it is obviously advantageous to study the both
subjects simultaneously in a unified geometric setting. Similar
constructions for other semi-Riemannian cases, complex hyperbolic
and quarternionic hyperbolic cases would also be very interesting.

{\noindent\bf Acknowledgement} The authors would like to thank to
Hyeonbae Kang, Chong Kyu  Han, Dohan Kim for help in the analysis.
In particular, the idea of $\epsilon$-approximation came across
with the suggestion of Professor Kang when we have a discussion
with him and then we formulate and develop as given in this paper.

\section{Hyperbolic sphere containing hyperbolic space}

Let $\Bbb R^{n,1}$ denote the real vector space $\Bbb R^{n+1}$ equipped with the bilinear form of signature $(n,1),$
$$
  \la x,y\ra =-x_0y_0+x_1y_1+ \cdots +x_ny_n,
$$
for all $x=(x_0,x_1,\cdots,x_n)$, $y=(y_0,y_1,\cdots,y_n)$.
Then the hyperbolic spaces $H^n_+$ and $H^n_-$, pseudo-sphere $S^n_1$ and light cone $L^n$ are defined by
$$
\al
 H^n_+&:= \{ x\in\Bbb R^{n,1} | \la x,x\ra =-1,\quad x_0>0\},\\
 H^n_-&:= \{ x\in\Bbb R^{n,1} | \la x,x\ra =-1,\quad x_0<0\},\\
 S^n_1&:= \{ x\in\Bbb R^{n,1} | \la x,x\ra =1\},\\
 L^n&:= \{ x\in\Bbb R^{n,1} | \la x,x\ra =0\}.
\eal
$$

We already know that $H^n_{\pm}$ is a Riemannian manifold which
has a constant sectional curvature $-1$, and that $S^n_1$ is a
Lorentzian  manifold (or semi-Riemannian of signature ($-,+,\cdots
,+$))  which has a constant  sectional curvature $1$ (see, for
example, \cite{O}). Now we consider  the Kleinian projective
model. By the radial projection $\pi_1$ with respect to the origin
from $H^n_+$ into $\{1\}\times \Bbb R^n$, we obtain the induced
Riemannian metric on the ball in $\{1\}\times \Bbb R^n$  as
follows (\cite{Geo},\cite{Ra}),
$$
ds^2_K=\left({\Sigma x_i dx_i\over 1-|x|^2}\right)^2 +{\Sigma dx_i^2\over 1-|x|^2}.
$$
We will denote the unit ball with this metric by $\Bbb H^n$.
If we extend this metric beyond  the unit ball using the same formula, this metric induces a semi-Riemannian structure on
the outside of the unit ball in $\{1\}\times \Bbb R^n$ and this will be denoted by $\Bbb L^n$. In fact, if we compare
this metric with the one induced from the Lorentzian space  $S^n_1\cap \{x\!=\!(x_0,x_1,\cdots ,x_n)|x_0\!>0\!\}$, by
the radial projection into  $\{1\}\times \Bbb R^n$, then they differ only by sign $-1$. This sign change of metric
implies the sign change of the sectional curvature  from $+1$ to $-1$, which, of course, the curvature of the metric
$ds^2_K$.
In this way,  we obtain an {\it extended Kleinian model} with metric $ds^2_K$ defined on $\{1\}\times \Bbb R^n$ except
 for the unit sphere $\{x\!=\!(x_1,\cdots ,x_n)||x|\!=\!(x_1^2+\cdots +x_n^2)^{1\over 2}\!=\!1\}$, and this extended
 hyperbolic space $( \{1\}\times \Bbb R^n, ds^2_K)$ will be denoted by $K^n$. Hence  $K^n$ consists of $\Bbb H^n,
 \Bbb L^n$, and the unit sphere.

Note that $K^n$ is simply an affine coordinate chart of the {\it projective space} $\Bbb R P^n_H$ equipped with the
singular metric composed of the usual hyperbolic metric coming from $H^n_+$ (or $H^n_-$) and the negative of the
Lorentz metric coming from $S^n_1$.

In this paper, it would be more convenient to consider the
Euclidean unit sphere in $\Bbb R^{n+1}$ with the induced
(singular) metric coming from this projective space  as a double
covering, and denote this model by  $\Bbb S^n_H$. This {\it
hyperbolic sphere model} $\Bbb S^n_H$ on the Euclidean sphere
$\{x\!=\!(x_0,x_1,$ $\cdots ,x_n)|x_0^2+x_1^2+\cdots +x_n^2\!=\!1\}$ has
three open  parts: Two radial images of $H^n_{\pm}$, called the
{\it hyperbolic part}, as two open disks on upper and lower
hemisphere and the radial image of $S^n_1$, called the {\it
Lorentzian part},   forms the remaining part which is
diffeomorphic to $\Bbb R^1\times \Bbb S^{n-1}$ (see Fig.2.1). Note
that the image of $S^1_1$ is composed of two components.

\begin{figure}[h]
\begin{center}
\includegraphics[width=0.7\textwidth]{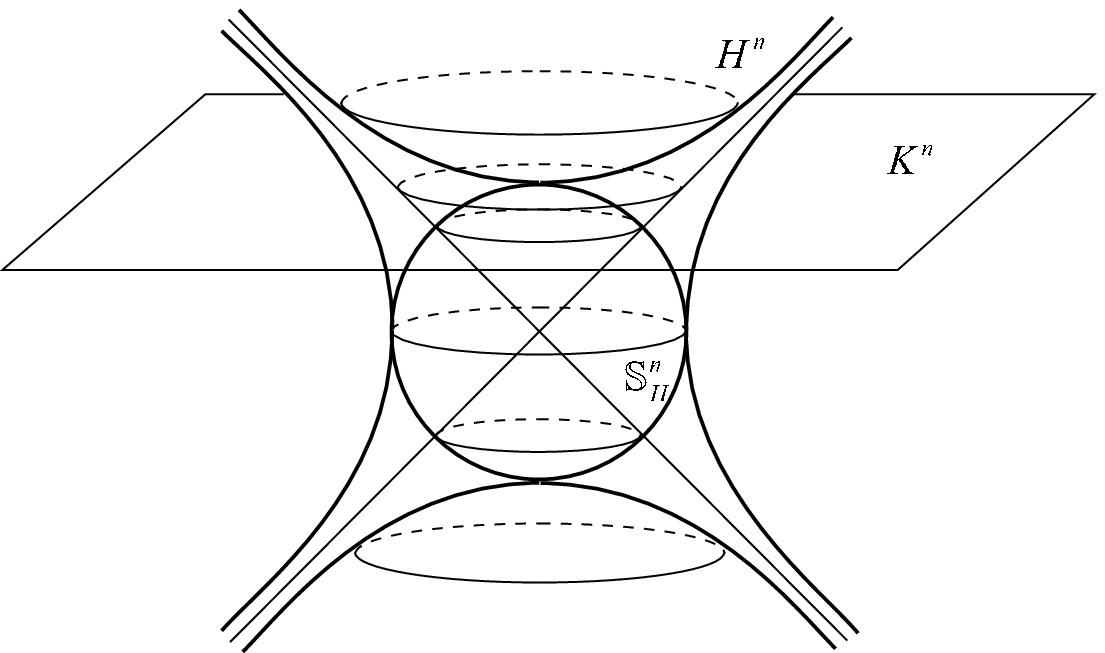}
\caption{\textbf{2.1}}  
\end{center}
\end{figure}

In this section,
we want to perceive this singular metric on $\Bbb S^n_H$ as an
analytic continuation of the hyperbolic part.
One way to do it is to view the metric as a limit of non-singular
complex metric.

We define an {\it $\epsilon-$approximation} of $ds^2_K$ as
$$
ds^2_{\epsilon}=\left({\Sigma x_i dx_i\over
d^2_{\epsilon}-|x|^2}\right)^2 +{\Sigma dx_i^2\over
d^2_{\epsilon}-|x|^2},
$$
where $d_{\epsilon}=1-\epsilon i$ with $\epsilon>0$ and
$i=\sqrt{-1}$. Then $ds^2_K=\lim_{\epsilon \to 0}ds^2_{\epsilon}$.

Since the length and angle can be viewed as an 1-dimensional
measure, we first study the volume (or measure) on $\Bbb S^n_H$
as a limit of the non-singular complex volume obtained from
$ds^2_{\epsilon}$. Note that the volume form on $K^n$ follows as
$$
\al
dV_K&=(\det(g_{ij}))^{\frac12}dx_1\wedge\cdots\wedge dx_n,\\
    &=\frac{dx_1\wedge\cdots\wedge dx_n}{(1-|x|^2)^{\frac{n+1}{2}}}.
\eal
$$
The volume form for $ds^2_{\epsilon}$ similarly will be given by
$$
dV_{\epsilon}=\frac{d_{\epsilon}dx_1\wedge\cdots\wedge dx_n}
{(d^2_{\epsilon}-|x|^2)^{\frac{n+1}{2}}}.
$$
In the Kleinian model, for a set $U$ in the unit disk $H^n$ the
volume of $U$ will be simply given by
$$
\text{vol} (U)=\int_U dV_K=\int_U  \frac{dx_1\wedge\cdots\wedge
dx_n}{(1-|x|^2)^{\frac{n+1}{2}}}
$$
as far as the integral exists. Furthermore this integral is also
obtained as $\lim_{\epsilon \to 0}\int_U dV_{\epsilon}$, since
$|d_{\epsilon}-|x|^2|\ge |1-|x|^2|$ and hence the Lebesgue
dominated convergence theorem applies. We have the same
conclusion for $U$ lying solely in the Lorentzian part $S^n_1$. But when we consider the volume function for Lorentzian
case, there is a sign problem.
There is a natural consistent way of choosing a sign $\{\pm i,\pm
1\}$ for the volume computation in the Lorentzian part once
$d_{\epsilon}=1-\epsilon i$ is chosen as our approximation, and
the sign convention will be explained later.

Now if a subset $U$ of $K^n$ lies across the light cone, that is
the boundary of $\Bbb H^n$, then the integral for $\text{vol} (U)$
does not make sense any more, and we want to define a volume of
$U$ as

\begin{equation}\label{x0} \mu(U)=\lim_{\epsilon \to 0}\int_U
\frac{d_{\epsilon}dx_1\wedge\cdots\wedge dx_n}
{(d^2_{\epsilon}-|x|^2)^{\frac{n+1}{2}}}
\end{equation}
whenever the limit exits. From now on we will call a Lebesgue
measurable set $U$ {\it $\mu$-measurable} if $\mu(U)$ is defined
as a finite value.

The actual computation of such integral doesn't seem to be easy.
But we can show that $\mu(U)$ is equal to the integral of the
singular volume form $dV_K$ over $U$ calculated in polar
coordinates if it is interpreted appropriately for a nice class
of the subsets $U$. One such case  is when the
integral can be considered as an analytic continuation in the
radial direction in the polar coordinates. Consider the integral for $\text{vol}
(U)$ as before.
$$
\al
\text{vol} (U)&=\int_U dV_K\\
              &=\int_U  \frac{dx_1\cdots dx_n}{(1-|x|^2)^{\frac{n+1}{2}}} \\
              &=\int_{G^{-1}(U)} \frac{r^{n-1}}{(1-r^2)^{\frac{n+1}{2}}} dr d\theta,
\eal
$$
where $G:(r,\theta)\mapsto(x_1,\cdots ,x_n)$ is the polar coordinates and  $d\theta$ is the volume form of the
Euclidean unit sphere $ \Bbb S^{n-1}$. If
$F(r)=\int_{G^{-1}(U)\cap S^{n-1}(r)}d\theta $ is an analytic function of $r$,
then we call such a subset $U$ of $K^n$ a {\it proper} set and the above integral becomes a 1-dimensional integral as
follows.
$$
 \int_{G^{-1}(U)} \frac{r^{n-1}}{(1-r^2)^{\frac{n+1}{2}}} dr d\theta=\int_a^b
\frac{r^{n-1}F(r)}{(1-r^2)^{\frac{n+1}{2}}} dr
$$
In general this integral does not make sense and the most natural thing we can do instead is to define $\text{vol}
(U)$ as the following contour integral
\begin{equation}\label{f2}
\text{vol} (U):=\int_{\gamma}
\frac{r^{n-1}F(r)}{(1-r^2)^{\frac{n+1}{2}}} dr
\end{equation}
where $\gamma$ is a contour from $a$ to $b$ for $a<1<b$ as depicted below. Here we will fix its contour type as
clockwise around $z=1$ once and for all throughout the paper.

\begin{figure}[h]
\begin{center}
\includegraphics[width=0.4\textwidth]{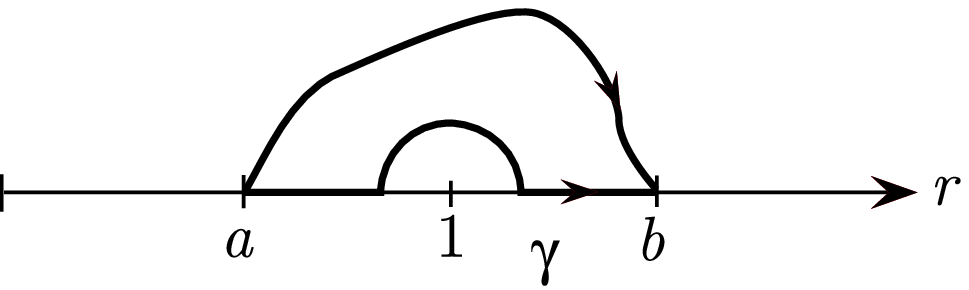}
\caption{\textbf{2.2}}
\end{center}
\end{figure}

Therefore we can compute the length of line segment on $\Bbb S^1_H$
by using the line integral (\ref{f2}). It is easy to see that
\begin{equation}\label{3}
m(b):=\int^b_{0,\gamma} \frac{d r}{1-r^2}=\left\{\alignedat2
&\frac {1}{2}\log \frac{1+b}{1-b},\qquad &0\le b<1,\\
&\frac {1}{2}\log \frac{b+1}{b-1} +\frac{\pi}{2}i,\quad &1<b,
                           \endalignedat\right.
\end{equation}
where the notation $\int^b_{a,\gamma}$ denote a line integral from
$a$ to $b$ along the contour $\gamma$.

Now we can show that this contour integral in fact gives us a way
of calculating $\mu(U)$.

\begin{pro}\label{1.1} For a proper set $U$ in $K^n$, $\mu(U)=\text{vol} (U)$. \
\end{pro}
\begin{proof} Using polar coordinates,
$$
\al \mu(U)&=\lim_{\epsilon \to 0}\int_{G^{-1}(U)}
\frac{d_{\epsilon}r^{n-1}}{(d_{\epsilon}^2-r^2)^{\frac{n+1}{2}}}
dr d\theta \\
&=\lim_{\epsilon \to 0}\int_a^b
\frac{d_{\epsilon}r^{n-1}F(r)}{(d_{\epsilon}^2-r^2)^{\frac{n+1}{2}}}
dr \eal
$$
for some $0\le a<1<b$. Then the equation $r^2-d_{\epsilon}^2=0$
gives poles $r=\pm d_{\epsilon}=\pm (1-\epsilon i)$, which, near
1, is located in the lower half of the complex plane $\Bbb C$.
Hence the integral $\int_a^b$ is equal to $\int_{\gamma}$ for any
contour $\gamma$ from $a$ to $b$ lying in the upper half of $\Bbb
C$. Now for such contour,
\begin{equation}\label{*}
\mu(U)=\lim_{\epsilon \to 0}\int_{\gamma}
\frac{d_{\epsilon}r^{n-1}F(r)}{(d_{\epsilon}^2-r^2)^{\frac{n+1}{2}}}
dr.
\end{equation}
In particular, choose $\gamma$ as depicted in Fig. 2.2 so that it
goes around $z=1$ through an upper semi-circle of radius $\delta$
in the clockwise direction. The upper semi-circle part of $\gamma$
can be given by
$$
\gamma(t)=1+\delta e^{(\pi-t)i}, \quad 0\le t \le \pi.
$$
Now for all sufficiently small $\epsilon < \epsilon_0$, choose
$\delta$ such that
$2\sqrt{\epsilon_0^4+4\epsilon_0^2}< \delta < 1$. Then we have
$$
\al |d_{\epsilon}^2-r^2|&=|-\epsilon^2-2\epsilon i
+1-r^2|\ge|1-r^2|-\sqrt{\epsilon^4+4\epsilon^2}\\
&\ge |1-r||1+r|-\frac{\delta}2 =  \delta|2+\delta
e^{(\pi-t)i}|-\frac{\delta}2 \ge
\delta-\frac{\delta}2=\frac{\delta}2. \eal
$$
This inequality gives us an uniform estimate of the integrand in
(\ref{*}) and completes the proof using Lebesgue dominated
convergence theorem.
\end{proof}

\begin{rem}\label{Remark} Our choice of $\epsilon$-approximation $d_{\epsilon}=1-\epsilon i$ is rather elaborate
anticipating proposition \ref{1.1}. If we choose $\tilde
d_{\epsilon}=1+\epsilon i$ instead, then $\tilde
\mu(U)=\lim_{\epsilon \to 0}\int_U d\tilde V_{\epsilon}$ will be
the same as the \~ {vol} $(U)$ which is given as in (\ref{x0})
where the contour type is given as {\it counterclockwise} around
$z=1$, i.e., lower going at $z=1$. The total length of $\Bbb
S^1_H$ is $2\pi i$ for the choice of $d_{\epsilon}$, and $-2\pi
i$ for the choice of $\tilde d_{\epsilon}$.
\end{rem}

To determine the various geometric quantities like the integration
for arc length on $\Bbb S^n_H$, the notion of norms of vectors is
essential. From the sign change of the  metric on the
pseudo-sphere $S^n_1$, the norm of  tangent vector
$x_p=(x_0,x_1,\ldots,x_n)$ at a Lorentzian part point is given by
$$\Vert x_p\Vert^2=-(-x_0^2+x_1^2+ \cdots +x_n^2),$$
and we should determine the sign of $\Vert x_p\Vert$ among $\pm 1$
and $\pm i$.

On the 1-dimensional hyperbolic sphere $\Bbb S^1_H$, we can
determine  the sign by looking at the sign of arc-length which can
be calculated by the formula (\ref{3}). Since the function $m(b)$
in the formula (\ref{3}) is decreasing on $b>1$,  the sign of
arc-length becomes negative value outside of $\Bbb H^1$ (see Fig.
2.3). This gives us ($-1$) as the right choice of the sign of
$\Vert x_p\Vert$ for the vectors in the radial direction at a
Lorentzian point.

For the sign for the vectors normal to the radial direction, we
need the following argument. It is not hard to check that the
clockwise contour integral of the volume form gives sign
$-i^{n-1}$ for Lorentzian part. Indeed for the volume integral in
the Lorentzian part will be given as
$$\int^b_a f(r) dr=\int^b_{0,\gamma} f(r) dr-\int^a_{0,\gamma} f(r) dr, \quad 1<a<b,$$
integrated along the clockwise contour $\gamma$. Here the
integrand $f(r)=\frac{r^{n-1}F(r)}{(1-r^2)^{\frac{n+1}{2}}}$ can
be written as $f(r)=\frac{g(r)}{(1-r^2)^{\frac{n+1}{2}}}$ with
$g(r)\ge 0$ for $r\ge 0$. The real axis for $r>1$ will have the
same sign under the map $f(r)$ and the sign can be traced using
the contour $\gamma$ which can be parametrized as $r=1+\delta
e^{(\pi -t)i}, 0\le t\le \pi$, near $r=1$. Then
$\frac{1}{(1-r^2)^{\frac{n+1}{2}}}=a e^{\frac{n+1}2 t i} (2+\delta
e^{(\pi -t)i})^{-\frac{n+1}{2}}, a=\delta^{-\frac{n+1}2}>0$, the
term $(2+\delta e^{(\pi -t)i})^{-\frac{n+1}{2}}$ doesn't give an
effect on the sign at $t=\pi$ and thus the sign at $t=\pi$ will be
$e^{\frac{n+1}2\pi i}=i^{n+1}=-i^{n-1}$. Hence on the
2-dimensional hyperbolic sphere $\Bbb S^2_H$, the volume for
Lorentzian part has the sign $-i$ and thus the consistent choice
of sign for the normal direction is $i$ (see Fig. 2.3). Now since
the dimension of the normal direction is $n-1$ with sign $i$ and
radial direction with sign $-1$, the volume has sign $(-1)\cdot
i^{n-1}=-i^{n-1}$, which  confirms again conceptually our volume
sign for the Lorentzian part.

It is easy to see that the light cone at a point in the Lorentzian
part is precisely the cone which is tangent to $\partial\Bbb H^n$
(see Fig. 2.3). Recall that our metric is the negative of the
metric of the pseudo-sphere $S^n_1$ and hence it has signature
$(+,-,\cdots ,-)$, which implies that $\Vert x_p\Vert^2>0$ inside
the light cone (the shaded region in Fig. 2.3) and $\Vert
x_p\Vert^2<0$ for normal directions. We can see this in the
Minkowski space by looking at the tangent vector $x_p$ at $p\in
S^n_1$ which lies inside and outside of the light cone at $p$
respectively.

\begin{figure}[h]
\begin{center}
\includegraphics[width=0.6\textwidth]{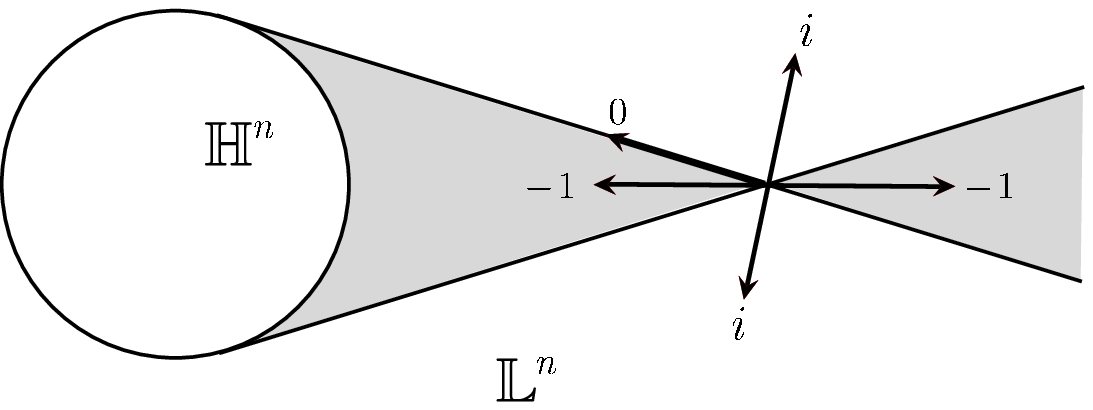}
\caption{\textbf{2.3}}
\end{center}
\end{figure}

We summarize the above discussion as the following convention.

\vskip 1pc \noindent{\bf Convention} {\it A tangent vector on the
hyperbolic part on $\Bbb S^n_H$ has a positive real norm, and a
tangent vector on the Lorentzian part on $\Bbb S^n_H$ has a
negative real, zero, or positive pure imaginary norm depending on
whether it is timelike, null, or spacelike respectively.} \vskip
1pc
 In \cite{Sc}, we can see a complex distance extending the cross
ratio to the exterior of the hyperbolic space. And we  can easily
check that the definition using  cross ratio and our definition
for distance and angle coincide. Our definition of angle will be
explained in Section 3.

We can see one of the similarities between the Euclidean sphere $\Bbb S^n$ and the hyperbolic sphere $\Bbb S^n_H$ in the
following theorem.
\vskip 1pc
\begin{thm}\label{1.4} vol($\Bbb S^n_H$)$=i^n\cdot$vol($\Bbb S^n$).
\end{thm}
\begin{proof} Putting $U=\Bbb S^n_H$ in the formula for
$\text{vol}(U)$ above, we have
\begin{equation}\label{3'}
\text{vol} (\Bbb S^n_H)=2(\int^{\infty}_{0,\gamma}
\frac{r^{n-1}}{(1-r^2)^{\frac{n+1}{2}}} dr) \text{ vol} (\Bbb
S^{n-1}).
\end{equation}
Similarly, we get the volume of the standard sphere by the radial projection from $\Bbb S^n$ to $\{1\}\times \Bbb R^n$, and
that is represented as
\begin{equation}\label{4}
\text{vol} (\Bbb S^n)=2(\int^{\infty}_0
\frac{t^{n-1}}{(1+t^2)^{\frac{n+1}{2}}} dt) \text{ vol} (\Bbb
S^{n-1}).
\end{equation}

\begin{figure}[h]
\begin{center}
\includegraphics[width=0.35\textwidth]{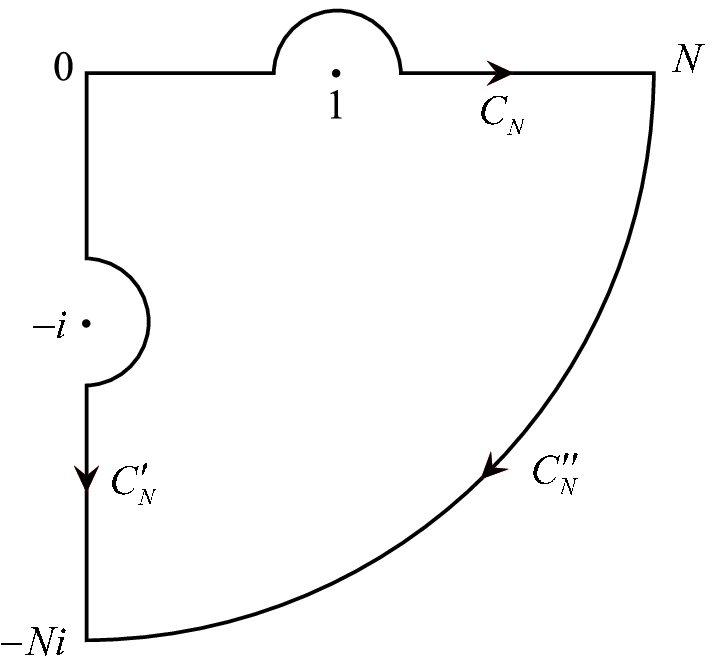}
\caption{\textbf{2.4}}
\end{center}
\end{figure}

To compare (\ref{3'}) and (\ref{4}), we need to define some
contours (see Fig. 2.4),

$$
\aligned
C_N:z_1(t)&=\left\{\aligned
                   &t,\qquad \qquad \quad 0\le t\le 1-\epsilon,\\
                   &1+\epsilon e^{\frac{i(1-t)\pi}{\epsilon}},  1-\epsilon\le t\le 1,\\
                   &t+\epsilon, \qquad\quad\ 1\le t\le N-\epsilon,\\
                   \endaligned \right.\\
C'_N:z_2(t)&=-iz_1(t),\\
C''_N:z_3(t)&=Ne^{-it},\qquad\qquad 0\le t \le \frac{\pi}{2}.\\
\endaligned
$$

For a given clockwise contour, we get
$$
\al
&\int^{\infty}_{0,\gamma} \frac{r^{n-1}}{(1-r^2)^{\frac{n+1}{2}}} dr\\
&=\lim_{N\to \infty} \int_{C_N} \frac{r^{n-1}}{(1-r^2)^{\frac{n+1}{2}}} dr \\
&=\lim_{N\to \infty} i^n \int_{C'_N}  \frac{t^{n-1}}{(1+t^2)^{\frac{n+1}{2}}} dt \qquad\quad r=it\\
&=\lim_{N\to \infty} i^n (\int_{C_N}  \frac{t^{n-1}}{(1+t^2)^{\frac{n+1}{2}}} dt +\int_{C''_N}
\frac{t^{n-1}}{(1+t^2)^{\frac{n+1}{2}}} dt)\\
&=\lim_{N\to \infty} i^n \int_{C_N}  \frac{t^{n-1}}{(1+t^2)^{\frac{n+1}{2}}} dt\\
&=i^n \int^{\infty}_0  \frac{t^{n-1}}{(1+t^2)^{\frac{n+1}{2}}} dt.
\eal
$$

From (\ref{3'}),(\ref{4}), and the above relation,  we deduce that
vol($\Bbb S^n_H$)$=i^n$ vol($\Bbb S^n$).
\end{proof}

\begin{rem}\label{1.5}
If we change the contour type of the integral (\ref{f2}), we have
different relation between vol($\Bbb S^n_H$) and vol($\Bbb S^n$).
If the contour is counterclockwise, in this case we should replace
$d_{\epsilon}$ to $\tilde d_{\epsilon}=1+\epsilon i$, then we have
vol($\Bbb S^n_H$)$=(-i)^n$ vol($\Bbb S^n$).
\end{rem}

\begin{rem}\label{1.6}
For various kinds of contour types, we easily deduce the following two formulas,
$$
\al
\text{vol }(\Bbb S^{2k-1}_H)&\equiv i^{2k-1} \text{vol }(\Bbb S^{2k-1}) \qquad (\text{mod }2 i^{2k-1}\text{ vol }
(\Bbb S^{2k-1})),\\
\text{vol }(\Bbb S^{2k}_H)&=i^{2k} \text{vol }(\Bbb S^{2k}).
\eal
$$
Above formulas say that the total volume of even dimensional model has unique value for any contour but odd dimensional
model has infinitely many values for various types of contours.
\end{rem}
For basic two contour types, we can go through into well established theory by using of $d_{\epsilon}$ and
$\tilde d_{\epsilon}$.

In the hyperbolic sphere $\Bbb S^n_H$, the hyperbolic part has $+\infty$ volume and the Lorentzian part has
$+\infty i^{n+1}$ ($+\infty (-i)^{n+1}$ for  $\tilde d_{\epsilon}$ case) volume. So the fact that the volume
of $\Bbb S^n_H$ is $i^n$ vol($\Bbb S^n$) gives the following nonsense equality.
$$\infty + \infty \cdot i^{n+1} = \text{vol }(\Bbb S^n)\cdot i^n $$ or
$$\infty + \infty \cdot (-i)^{n+1} = \text{vol }(\Bbb S^n)\cdot (-i)^n, \text{  for  }  \tilde d_{\epsilon} \text{  case}.$$
But the nonsense equality does not come true in a suitable measure
theory, hence we need not worry about that. In Section 3, we can
see the suitable and nice finitely additive complex measure theory
for the spaces $K^n, \Bbb R P^n_H$, and $\Bbb S^n_H$.

\section{Invariant property and measure theory on the extended
model}

Let us study the invariance of $\mu$-measurable set in $\Bbb R
P^n_H$ or $\Bbb S^n_H$. Here the isometry group of $\Bbb S^n_H$ is
$O(n,1)$ and we identify $PO(n,1)$ as the subgroup of $O(n,1)$
leaving $H^n_+$ invariant. If $U$ is a measurable set in the
hyperbolic part or in the Lorentzian part, then the $\mu$-measure is
simply the usual hyperbolic or Lorentzian volume (with appropriate
sign) and hence the invariance is obvious. But if we consider $U$
lying across the $\partial \Bbb H^n$, the boundary of the
hyperbolic part, then the invariance seems to be a subtle problem.
Even if we do not know the invariance for general $\mu$-measurable
set holds, we are able to show the invariance for a region with
piecewise analytic boundary transversal to $\partial\Bbb H^n$.
More precisely, we define as follows.

\begin{defi}\label{2.0}
An $n$-dimensional region $U$ of the extended hyperbolic space
$\Bbb S^n_H$ is called a region with piecewise analytic boundary
transversal to $\partial\Bbb H^n$ if it is bounded by finitely
many analytic sets, $Ai:=\{x|f_i(x)=0, f_i \text{ is an analytic
function }\}$, $i=1,2,\ldots,k$, which are transversal to
$\partial\Bbb H^n$ such that for each $p\in U\cap\partial \Bbb
H^n$, $p$ belongs to one of the bounding analytic sets but at most two
such sets $A_i$ and $A_{i+1}$ (mod $k$) with the property that
$\nabla f_i(p)$, $\nabla f_{i+1}(p)$ and the normal vector to
$\partial\Bbb H^n$ are linearly independent.
\end{defi}

\begin{pro}\label{2.1} Let $U$ be a region with piecewise analytic boundary transversal to $\partial\Bbb H^n$ in the
extended hyperbolic space. Then $\mu(U)$ has a well-defined finite
value, i.e., $U$ is $\mu$-measurable, and $\mu (g(U))=\mu (U)$ for
each $g\in PO(n,1)$.
\end{pro}

\begin{figure}[h]
\begin{center}
\includegraphics[width=0.6\textwidth]{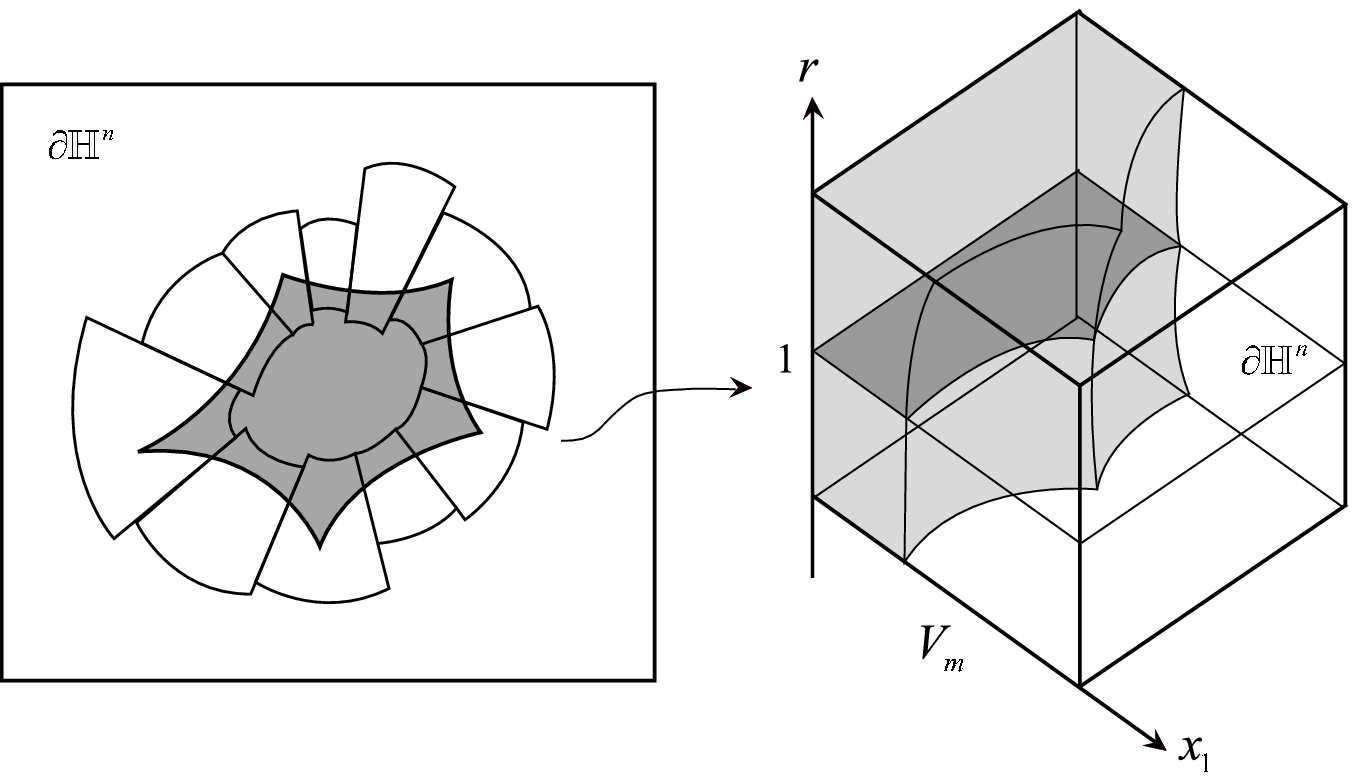}
\caption{\textbf{3.1}}
\end{center}
\end{figure}

\begin{proof}We want to show first that $F_U(r)=\int_{G^{-1}(U)\cap
S^{n-1}(r)}~d\theta$ and $F_{gU}(r)$ are real analytic function of
$r$ near $r=1$, i.e., $U$ and $gU$ are proper sets, (see the
previous section for the notations and definition) so that $U$ and
$gU$ have well-defined finite values.

We decompose a neighborhood of $U$ near $r=1$ appropriately so that
each piece has a box coordinate under a real analytic coordinate
change as in Figure 3.1. If we denote two analytic boundaries in
$V_m$ by $f_m(x_1,\ldots,x_{n-1},r)=0$ and
$f_{m+1}(x_1,\ldots,x_{n-1},r)=0$, then we can rewrite these as
$x_1=g_m(x_2,\ldots,x_{n-1},r)$ and
$x_1=g_{m+1}(x_2,\ldots,x_{n-1},r)$ by the (analytic) implicit
function theorem. From the condition that $\nabla
f_m(x_1,\ldots,x_{n-1},1)$ and $\nabla
f_{m+1}(x_1,\ldots,x_{n-1},1)$ are not parallel, we may assume
$\frac{\partial}{\partial x_2} g_m\ne\frac{\partial }{\partial
x_2}g_{m+1}$, and then again by the (analytic) implicit function
theorem, we can write the intersection of two analytic sets,
$g_m(x_2,\ldots,x_{n-1},r)=g_{m+1}(x_2,\ldots,x_{n-1},r)$, as
$x_2=h(x_3,\ldots,x_{n-1},r)$.

Now the function $F_U(r)$ for the shaded region $U$ in $V_m$ in
Fig. 3.1 will be given by
$$\al
F_U(r)&=\int^{b_{n-1}}_{a_{n-1}}\cdots\int^{b_3}_{a_3}\int^{h(x_3,\ldots,x_{n-1},r)}_{a_2}\int^{g_m(x_2,\ldots,
x_{n-1},r)}_{a_1}~f(x)~dx_1 dx_2 \cdots dx_{n-1}\\
&+
\int^{b_{n-1}}_{a_{n-1}}\cdots\int^{b_3}_{a_3}\int_{h(x_3,\ldots,x_{n-1},r)}^{b_2}\int^{g_{m+1}(x_2,\ldots,
x_{n-1},r)}_{a_1}~f(x)~dx_1 dx_2 \cdots dx_{n-1}.\eal$$
 Since the function $h,g_m,g_{m+1}$, and $f$ are analytic,
 $F_U(r)$ is analytic as desired. Similarly $F_{gU}(r)$ is also
 analytic.

Next we show that $\mu (U)=\mu (gU)$. Let $\gamma$ be a contour
from $a$ to $b$ in Fig. 2.2 and let $B_r=G^{-1}(U)\cap
S^{n-1}(r)$. Then
$$\al
\mu(U)&=\lim_{\epsilon\to 0} \int_U
dV_{\epsilon}=\lim_{\epsilon\to 0} \int^b_a
\frac{r^{n-1}F_U(r)}{(d^2_{\epsilon}-r^2)^{\frac{n+1}{2}}} dr\\
&=\lim_{\epsilon\to 0} \int_{\gamma}
\frac{r^{n-1}F_U(r)}{(d^2_{\epsilon}-r^2)^{\frac{n+1}{2}}}
dr=\int_{\gamma} \frac{r^{n-1}F_U(r)}{(1-r^2)^{\frac{n+1}{2}}}
dr=\int_{\gamma}\int_{B_r} dV_0. \eal
$$
Here the third equality holds since the pole of the integrand near
1 has a negative imaginary part. The fourth equality follows from
$|d^2_{\epsilon}-r^2|\ge \frac{\delta}{2}$ as in the proof of
Proposition 2.1 so that the integrand is uniformly bounded for all
sufficiently small $\epsilon>0$ and we can apply Lebesgue
dominated convergence theorem.

Now if we can show for $\mu (gU)$ similarly as follows, then the
proof will be completed.
$$\al
\mu(gU)&=\lim_{\epsilon\to 0} \int_{gU}
dV_{\epsilon}=\lim_{\epsilon\to 0} \int_{U}
g^*(dV_{\epsilon})=\lim_{\epsilon\to 0} \int^b_a \int_{B_r}
g^*(dV_{\epsilon})\\&=\lim_{\epsilon\to 0} \int_{\gamma}
\int_{B_r} g^*(dV_{\epsilon})=\int_{\gamma} \int_{B_r}
g^*(dV_0)=\int_{\gamma} \int_{B_r} dV_0.
 \eal
$$
Here $g$ is an isometry and hence $g^*(dV_0)=dV_0$, but notice
that $g$ does not preserve $\epsilon$-volume form $dV_{\epsilon}$
and we need several steps as before. We have to check 4th and 5th
equality and all others are obvious.

For $g\in PO(n,1)$, if we let
$dV_{\epsilon}=f_{\epsilon}(r,\theta)~d\theta dr$, then
$$g^*(dV_{\epsilon})=(f_{\epsilon}\circ g)(\det g')~d\theta dr=:h_{\epsilon}(r,\theta)~d\theta dr.$$
Since $g$ is an isometry, $g$ obviously preserves the orientation
of real $r$-axis. When we complexify everything, as an analytic
map $g$ preserves the orientation of complex $r$-axis and leaves
real $r$-axis invariant being a real map. Hence the pole of
$h_{\epsilon}$ also lies below the real axis as $f_{\epsilon}$
does.

Now for the 5th equality, note that $dV_{\epsilon}$ and
$g^*(dV_{\epsilon})$ differ only by Jacobian determinant for their
values. Since we already know that $dV_{\epsilon}$ is uniformly
bounded on its domain of integral $\gamma\times B_r$  so is
$g^*(dV_{\epsilon})$ and we can apply Lebesgue dominated
convergence theorem to obtain the equality.
\end{proof}

\begin{rem}\label{2.2} In fact, we can generalize this proposition to a region with more
general type boundaries than analytic boundaries. (See \cite{CK}.)
\end{rem}

Let $\mathcal F$ be the collection of $\mu$-measurable subsets of
$K^n$. Then it is obvious that $\mu: \mathcal F \to \Bbb C$ is a
finitely additive set function, i.e., for $E_1,E_2\in \mathcal F$
and disjoint union $E_1 \amalg  E_2$, we have $\mu(E_1 \amalg
E_2)= \mu(E_1)+\mu(E_2)$ so $E_1 \amalg  E_2 \in \mathcal F$ (this
result trivially generalizes to that if two sets among $E_1,E_2,E_1
\amalg  E_2$ are $\mu$-measurable, then the other is
$\mu$-measurable).

It would be nice if $\mu: \mathcal F \to \Bbb C$ is a finitely
additive measure on an algebra $\mathcal F$. But unfortunately
$\mathcal F$ is not an algebra, as we see in the following
example.

 Consider a cone $E_1$ over a small disk $B$ on $\partial \Bbb H^n$ with the vertex at the origin of
 $\Bbb H^n$, where the origin is the point $\{1\}\times \{0\}\in \{1\}\times \Bbb R^n=K^n$. Let $g$ be a rotation
 such that $g(E_1\cap \Bbb H^n)$ is disjoint from $E_1\cap \Bbb H^n$, and let $E_2=g(E_1\cap \Bbb H^n)\cup
 (E_1\cap \Bbb L^n)$. Then $E_1$ and $E_2$ are clearly $\mu$-measurable sets and by Theorem \ref{1.4}, for $i=1,2$,
$$
\al
\mu(E_i)=\text{vol }(E_i)&=\left( \int_{0,\gamma}^{\infty} \frac{r^{n-1}}{(1-r^2)^{\frac{n+1}{2}}}dr\right) \text{vol }(B)\\
&=\frac{i^n}{2} \frac{\text{vol }(\Bbb S^n)}{\text{vol }(\Bbb S^{n-1})}\text{vol }(B)=\frac{i^n\sqrt{\pi}}{2}
\frac{\Gamma(\frac n 2)}{\Gamma(\frac {n+1}2)}\text{vol }(B).
\eal
$$
But $E_1\cap E_2$ and $E_1\cup E_2$ have infinite volume, hence
not $\mu$-measurable. This example looks rather artificial, but we
will show that the problem of finding a useful large enough
algebra in $\mathcal F$ is a very sensitive  problem in its
nature.

The simplest type of $\mu$-measurable sets in $\Bbb H^n$ (or in
$\Bbb L^n$) is bounded regions not touching the $\partial \Bbb
H^n$, because the  $\mu$-measurable sets in $\Bbb H^n$ (or in
$\Bbb L^n$) is exactly the same as the original positive volume
measure in $\Bbb H^n$ (or in $\Bbb L^n$ up to sign).  Also the
cone region with vertex at the origin of $\Bbb H^n$ is
$\mu$-measurable.  Our first attempt to finding an algebra in
$\mathcal F$  is the one generated by sets of the above three
types.

Let $\mathcal U_h$ (resp. $\mathcal U_l$) be the collection of
$\mu$-measurable sets in  $\Bbb H^n$ (resp. $\Bbb L^n$) and
$\mathcal U_c$ be the collection of $\mu$-measurable cones with
vertex at the origin of $\Bbb H^n$. We assume that $\varnothing\in
\mathcal U_c$. Let $\mathcal M$ be the smallest algebra containing
$\mathcal U_h, \mathcal U_l$ and $\mathcal U_c$, so that $\mathcal
M$ contains $\mathcal U_d:=\{V|~V=U_2-U_1, U_2\in \mathcal U_c
\text{ and  } U_1=U_h\cup U_l, \text{ where }U_h\in\mathcal U_h
\text{ and } U_l\in \mathcal U_l\}.$

\begin{pro}\label{2.3} (1) Let $\mathcal M'=\{V_1\cup V_2|~ V_2\in \mathcal U_d \text{  and  }
 V_1=V_h\cup V_l, \text{ where }V_h\in\mathcal U_h \text{ and }
V_l\in \mathcal U_l
 \}$ and $\mathcal M''$ be the subcollection of
$\mathcal M'$ with the condition $V_1\cap V_2=\varnothing$. Then we
have $\mathcal M=\mathcal M'=\mathcal M''$.

\noindent(2) For any $U\in\mathcal M$, $U$ can be written as a
disjoint union $U=U_h\amalg U_l \amalg U_d$, where $U_h\in
\mathcal U_h$, $U_l\in \mathcal U_l$ and $U_d\in \mathcal U_d$.
Hence the collection $\mathcal M$ is a $\mu$-measurable algebra.
\end{pro}

\begin{proof} An element of $\mathcal M'$ can be written as
$V_1\cup(U_2-U_1)$ and the proof that $\mathcal M'=\mathcal M''$
follows readily from the  general set identity:
$$V_1 \cup(U_2-U_1)=(V_1-U_2)\amalg (U_2-(U_1-V_1)).$$

We want to show that $\mathcal M=\mathcal M'$. Since $\mathcal
M'\subset \mathcal M$ and $\mathcal U_h, \mathcal U_l, \mathcal
U_c\subset \mathcal M'$, it suffices to show that $\mathcal M'$ is
an algebra. The whole space $K^n\in\mathcal M'$ since it is a
cone. It is easy to show that $V\in\mathcal M'$ implies
$V^c\in\mathcal M'$. Indeed this follows from the general set
identity:
$$(V_1 \cup(U_2-U_1))^c=(U_1-V_1)\cup (U_2^c-V_1).$$
Note that the complement of a cone is also a cone.

Finally, let's show that $\mathcal M'$ is closed under the union
operation. Write
$$(V_1 \cup(U_2-U_1))\cup(V'_1 \cup(U'_2-U'_1))=(V_1\cup V'_1)\cup ((U_2-U_1)\cup(U'_2-U'_1)).$$
and observe the set identity:
$$(U_2-U_1)\cup(U'_2-U'_1)=(U_2\cup U'_2)-((U_1- U'_2)\cup (U'_1-U_2)\cup(U_1\cap U'_1)).$$

Now (2) follows from (1) trivially.
\end{proof}

\begin{figure}[h]
\begin{center}
\includegraphics[width=0.4\textwidth]{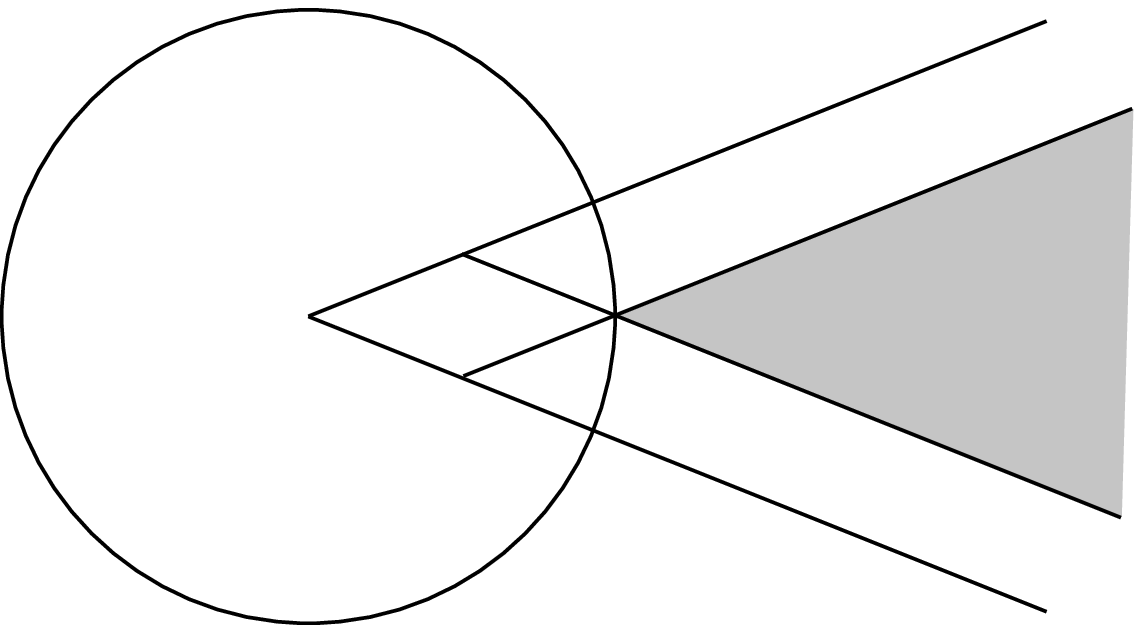}
\caption{\textbf{3.2}}
\end{center}
\end{figure}

In dimension 2, $\mathcal M$ gives us a fairly large class of
$\mu$-measurable sets. We know that a cusp in $\Bbb H^2$ has a
finite area. The same thing holds for a cusp in $\Bbb L^2$. This
can be shown easily by a direct computation of the integral in the
plane $K^2$, or by looking at a configuration of three cones
depicted in Fig. 3.2, if we notice that any cone with vertex in
$\Bbb H^2$ has a finite area by proposition \ref{2.1}.

These observations and that $\mu$-measure is essentially a positive
measure in $\Bbb H^n$ and $\Bbb L^n$ lead us immediately to
conclude that a domain transversal to $\partial\Bbb H^2$ belongs
to $\mathcal M$ without using Proposition \ref{2.1} (see Fig.
3.3).

\begin{figure}[h]
\begin{center}
\includegraphics[width=0.5\textwidth]{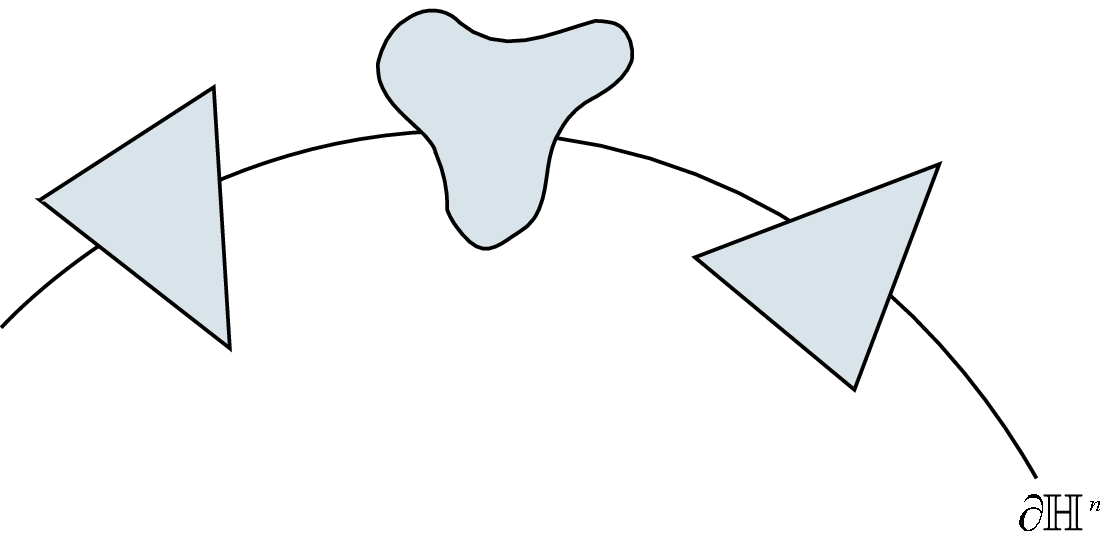}
\caption{\textbf{3.3}}
\end{center}
\end{figure}

But if the dimension is greater than or equal to 3, the algebra
$\mathcal M$ becomes a small collection which doesn't even contain
the class of cones with vertex in $\Bbb H^n$, whose volume, as we
already know, are finite.

Indeed let's consider the following two cones in $\Bbb R^n=K^3$.
$$
\al
C_1: z&\ge\frac3 4 r\\
C_2: z&\ge\frac1 3 (r+1)
\eal
$$

\begin{figure}[h]
\begin{center}
\includegraphics[width=0.5\textwidth]{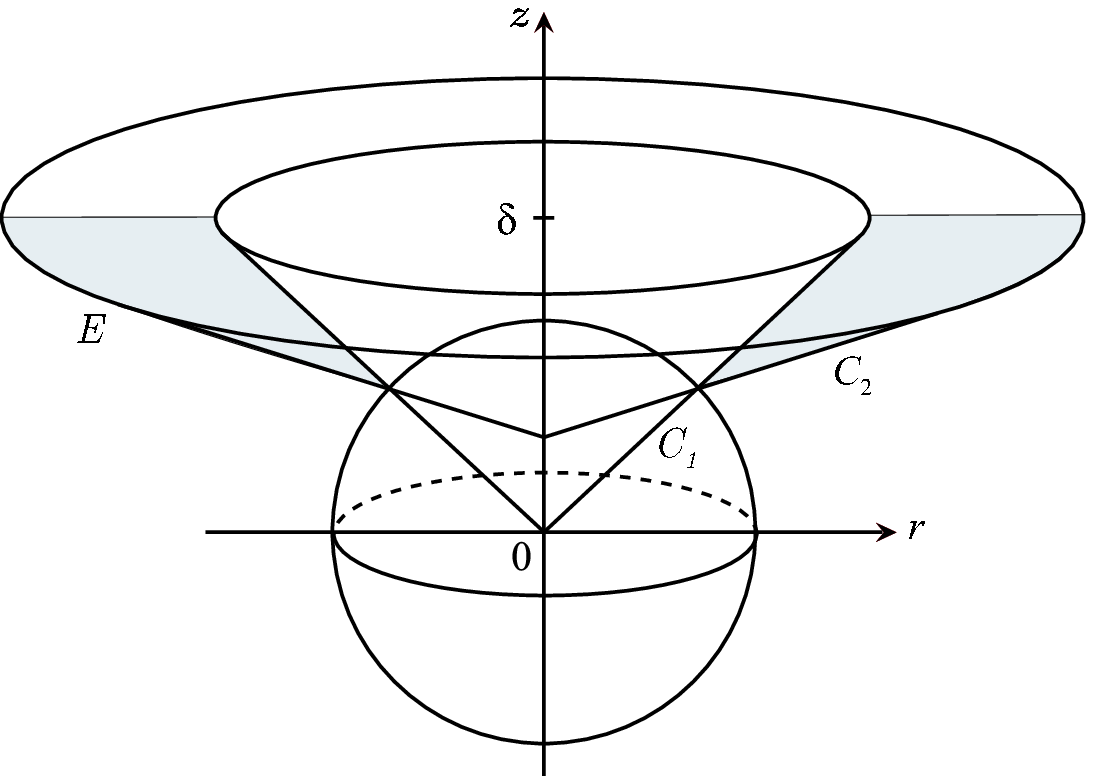}
\caption{\textbf{3.4}}
\end{center}
\end{figure}

Then $C_1\in \mathcal M$ by the definition of $\mathcal M$, but we
claim $C_2\notin \mathcal M$ by showing that the volume of
$E=C_2-C_1$ is infinite by direct computation.
$$
\al
\mu (E)&=\iiint_E \frac {dx dy dz}{(1-x^2-y^2-z^2)^2}\\
&=\int_{\frac 3 5}^{\delta}\int_{0}^{2\pi}\int_{\frac 4 3 z}^{3z-1}\frac {r dr d\theta dz}{(1-z^2-r^2)^2}\\
&=\pi\int_{\frac 3 5}^{\delta}\frac {1}{z-\frac 3 5}\left(\frac{9}{25}\frac{1}{z+\frac 3 5}-\frac{1}{10 z}\right) dz.
\eal
$$
This integral  diverges since the function in parenthesis in the integrand has a positive lower bound on
$[\frac 3 5 ,\delta]$.

Even though the collection $\mathcal M$ does not contain natural
domains we want in dimensions greater than or equal to 3,  we can construct
sufficiently large and very important collections $\mathcal H$ and
$\mathcal H'$ in such dimensions.

\begin{pro}\label{2.4} Let $\mathcal H$ be the algebra in $K^n$ (resp. $\Bbb S^n_H$) generated by half spaces (resp.
 hemisphere) not tangent to $\partial \Bbb H^n$. Then the
collection $\mathcal H$ is a $\mu$-measurable algebra.
\end{pro}

\begin{proof}
To show $\mu$-measurability of $U\in\mathcal H$, it suffices to
consider only $\mu$-measurability of the region
$U\cap\{1-\delta\le r\le 1+\delta\}$ for a small $\delta>0$ since
the remaining parts are completely contained in either $\Bbb H^n$
or $\Bbb L^n$.

\begin{figure}[h]
\begin{center}
\includegraphics[width=0.9\textwidth]{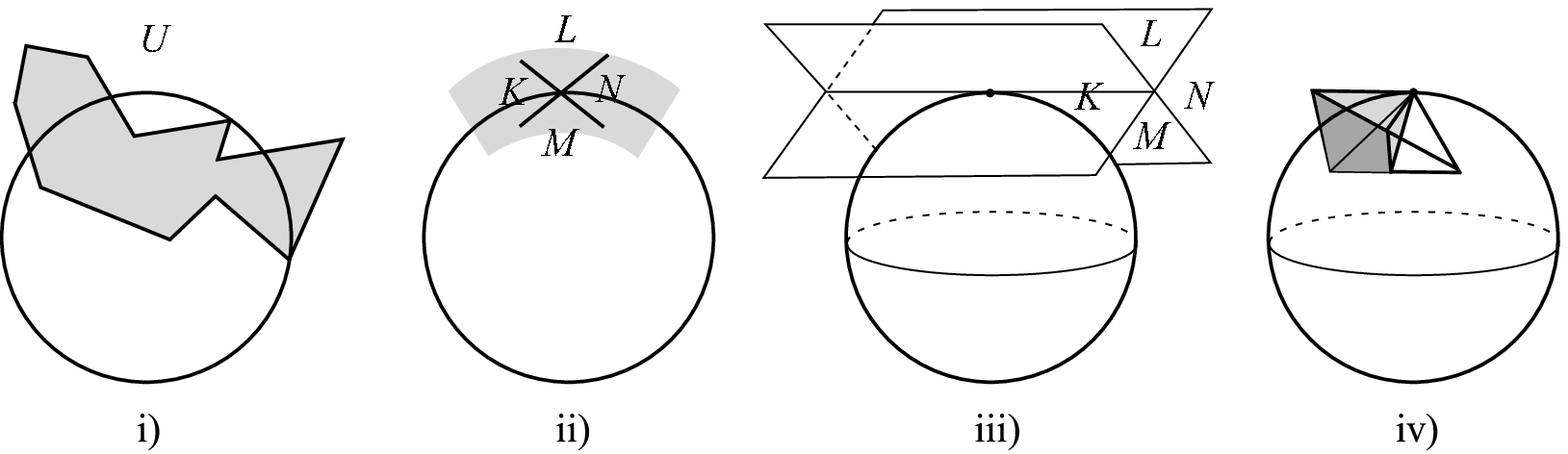}
\caption{\textbf{3.5}}
\end{center}
\end{figure}

When the dimension $n=2$, by virtue of Proposition 3.2, the only
exceptional case  we have to worry are the cone types regions
$K,L,M,N$ given as in Fig 3.5-(ii). We know that the hyperbolic
cusp $M$ and the Lorentzian cusp $L$ have finite areas.

The key observation about the $\mu$-measurability we are using
over and over in the following argument is that if two of $A,B$
and $A\cup B$ are $\mu$-measurable, so is the remaining one. For
instance, $K$ is $\mu$-measurable since $L$ and $K\cup L$ are
$\mu$-measurable. In this way all these four types of region are
$\mu$-measurable.

When $n=3$, we have to consider several exceptional case.

Case 1. When two planes transversal to $\partial \Bbb H^3$ meet at
a line which is tangent to $\partial \Bbb H^3$: (See Fig.
3.5-(iii).) Again in this case, we have four regions $K,L,M,N$ as
before and we can show that the region $L$ contained in the
Lorentzian part has a finite volume by direct computations (see
APPENDIX). Then the other regions have finite volumes by the same
argument as the case $n=2$ using the key observation.

Case 2. When three transversal planes meet at a point on $\partial
\Bbb H^3$: In this case, we obtain 8 regions from these 3 planes.
At least one of these regions is contained in the Lorentzian part
and this is contained in an $L$-type region considered in Case 1.
This implies that it has a finite volume since $\mu$-measure on
the Lorentzian part is essentially a positive measure. Then a
neighboring octant has a finite volume by the key observation and
Proposition 3.2. Hence again their neighboring octants have finite
volumes by the same reasoning and so forth.

Case 3. When four or more transversal planes meet at a point on
$\partial \Bbb H^3$: Fig. 3.5-(iv) shows a region (shaded one)
bounded by 4 transversal planes. This region is the difference of
two regions of the type considered in Case 2 and hence is of
finite volume again by the key observation. Now finiteness of
volume of a region bounded by many planes follows by induction.

If $n\ge 4$, we have more exceptional cases but still we have the
same conclusion by the same argument and observations since  we
can start with a piece in the Lorentzian part which has a finite
volume by APPENDIX.
\end{proof}

\begin{pro}\label{2.5} Let $\mathcal H'$ be the smallest algebra in $K^n$ or $\Bbb S^n_H$ containing
$\mathcal U_h, \mathcal U_l$, and $\mathcal H$. Then the
collection $\mathcal H'$ is $\mu$-measurable algebra.
\end{pro}

\begin{proof} Use the same method shown in Proposition \ref{2.3}.
\end{proof}

In dimension 2, the collection $\mathcal M$ is strictly larger
than $\mathcal H'$, i.e., $\mathcal H' \subsetneq \mathcal M$. We
will leave the proof that $\mathcal H' \subsetneq \mathcal M$ as an
easy exercise for readers.

\section{Lengths and angles on the extended
hyperbolic space}

We denote  the distance between two points $A$ and $B$ in the
extended hyperbolic space $\Bbb S^n_H$ as $d_H(A,B)$. Let's first
discuss the distance between two points on $\Bbb S^1_H$. In this
case, the $\epsilon$-metric and the volume form are given as
$ds_{\epsilon}^2=\frac{d_{\epsilon}^2 dx_1^2}
{(d^2_{\epsilon}-x_1^2)^2}$ and $dV_{\epsilon}=\frac{d_{\epsilon}
dx_1} {d^2_{\epsilon}-x_1^2}=ds_{\epsilon}$ on its affine chart
$K^1$. If we let the affine coordinates of two points $A$ and $B$
of $\Bbb S^1_H$ be $x_1=a$ and $x_1=b$ respectively, then the
$\mu$-measure of the line segment $l=[a,b]$ is
$$\mu(l)=\lim_{\epsilon \to 0}\int_l ds_{\epsilon}=\lim_{\epsilon \to
0}\int^b_a  \frac{d_{\epsilon} dx_1} {d^2_{\epsilon}-x_1^2},$$ and
this is to be the distance $d_H(A,B)$ of $A$ and $B$. For
instance, if $A$ and $B$ are symmetric with respect to the light
cone $x_0=x_1$ in $\Bbb R^{1,1}$ as in Fig. 4.1. (i.e., $A$ and
$B$ as vectors of $\Bbb R^{1,1}$ are perpendicular), then their
affine coordinates are $a(<1)$ and $\frac 1a$, and the distance
will be $\frac{\pi}2i$ by formula (\ref{3}). The distance between
isometric images $A'$ and $B'$ of $A$ and $B$ will be again
$\frac{\pi}2i$ being symmetric, and hence $d_H(B,B')=-d_H(A,A')$
in Fig. 4.1.

\begin{figure}[h]
\begin{center}
\includegraphics[width=0.35\textwidth]{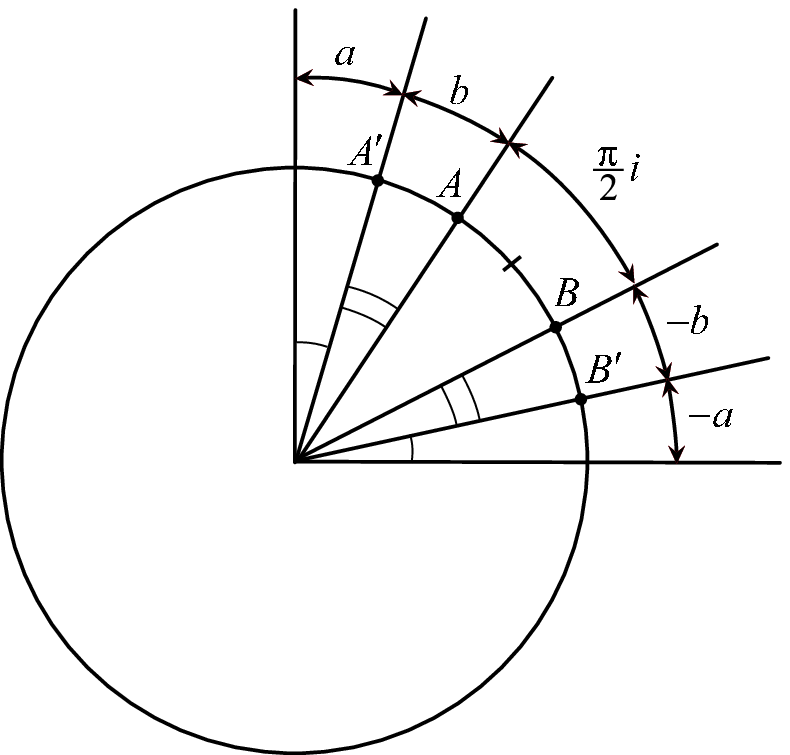}
\caption{\textbf{4.1}}
\end{center}
\end{figure}

To discuss the distance between two points in $\Bbb S^n_H$ in
general, it suffices to consider on $\Bbb S^2_H$. It is natural to
define a distance between two points on $\Bbb S^2_H$ by the
$\mu$-measure of the geodesic line segment connecting these two
points.

For actual computations, it would be convenient to divide into the
following 3 cases. For the case when the geodesic connecting two
points meet $\partial\Bbb H^2$ transversely, we may assume that
these two points lie on $\Bbb S^1_H=\Bbb S^2_H\cap \{x|x_2=0\}$ by
an isometry and can handle as discussed above.

For the case when the geodesic line connecting these two points
does not intersect $\partial\Bbb H^2$, we can send this line to
the equator $(=\Bbb S^2_H\cap \{x|x_0=0\})$ of $\Bbb S^2_H$ by an
isometry, and hence the distance becomes $i$ times the distance on
the standard Euclidean unit circle.

The remaining case is when the line is tangent to $\partial\Bbb
H^2$. On $K^2$, if we restrict the Kleinian $\epsilon$-metric to
the line given by $x_2=k$, then we have
$ds_{\epsilon}^2=\frac{(d_{\epsilon}^2-k^2) dx_1^2}
{(d^2_{\epsilon}-k^2-x_1^2)^2}$ and the length of the line segment
connecting $(0,k)$ and $(a,k)$ for $a>0$ will be given by
$\lim_{\epsilon \to 0}\int_0^a \frac{(d_{\epsilon}^2-k^2)^{1\over
2} dx_1} {d^2_{\epsilon}-k^2-x_1^2}$.  Since the length of the
line segment in $\Bbb H^2$ is positive, we have to choose one
whose real part is positive among the two possible values of
$(d_{\epsilon}^2-k^2)^{1\over 2}$ when $0<k<1$. By continuity we
choose $c_{\epsilon}=(d_{\epsilon}^2-1)^{1\over 2}$ with
$Re((d_{\epsilon}^2-1)^{1\over 2})>0$ (or
$Im((d_{\epsilon}^2-1)^{1\over 2})<0$). Now the length of the line
segment connecting $(0,1)$ and $(a,1)$ for $a>0$ will be
$$\al&\lim_{\epsilon \to 0}\int_0^a \frac{(d_{\epsilon}^2-1)^{1\over
2} dx_1} {d^2_{\epsilon}-1-x_1^2}=\lim_{\epsilon \to 0}\int_0^a
\frac{c_{\epsilon} dx_1} {c^2_{\epsilon}-x_1^2}=\lim_{\epsilon \to
0}\frac 12 \log\frac{c_{\epsilon}+x_1} {c_{\epsilon}-x_1}\Big
|^a_0\\
=&\lim_{\epsilon \to 0}\frac 12 \log\frac{c_{\epsilon}+a}
{c_{\epsilon}-a}=\lim_{\epsilon \to 0}\frac 12
\log\frac{r_1e^{i\theta_1}} {r_2e^{i\theta_2}}=\lim_{\epsilon \to
0}(\frac 12 \log\frac{r_1}
{r_2}+\frac{i(\theta_1-\theta_2)}{2})=\frac{\pi}{2}i \eal$$

Now from the finite additivity of $\mu$-measure, we conclude as
follows.

\begin{lem}\label{3.1}For a point $x$ lying on $\partial \Bbb H^2$ and a dual geodesic $x^{\bot}$, the lengths of the
line segments in $x^{\bot}$ are defined by
$$
\al
d_H(w,y)&=0, \quad\text{if } w,y \text{ are in the same side with respect to }x,\\
d_H(y,z)&=\pi i, \quad\text{if } y,z \text{ are in the opposite sides with respect to }x,\\
d_H(x,y)&=d_H(x,z)=\frac{\pi}{2}i. \eal
$$
See the following Fig. 4.2.
\end{lem}

\begin{figure}[h]
\begin{center}
\includegraphics[width=0.5\textwidth]{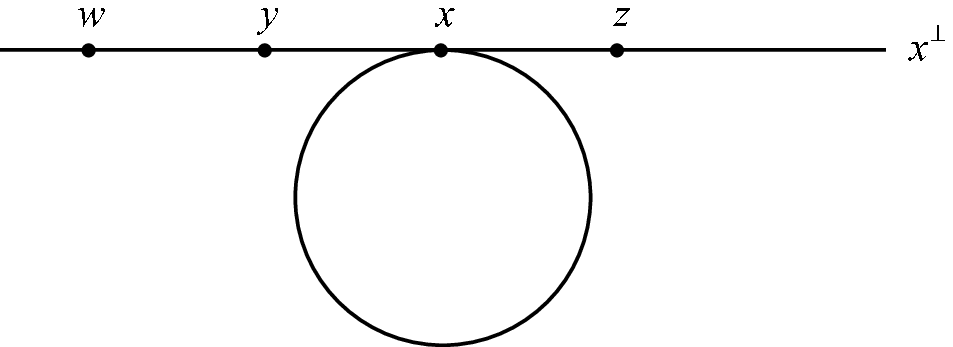}
\caption{\textbf{4.2}}
\end{center}
\end{figure}

In all of these discussions, we in fact have to show that
1-dimensional $\mu$-measure is invariant under isometry since the
hyperbolic isometry does not preserve $\epsilon$-metric. We will
show this in the following Theorem \ref{3.3} below.

We conclude the following theorem from the above discussion.
\begin{thm}\label{3.2} The total length, i.e., 1-dimensional $\mu$-measure, of any great circle in $\Bbb S^n_H$
is $2 \pi i$.
\end{thm}

The intersection of $(k+1)$-dimensional subspace of $\Bbb R^{n,1}$
with $\Bbb S^n_H$ is a totally geodesic $k$-dimensional subspace
of $\Bbb S^n_H$ and the $\epsilon$-metric on this $k$-dimensional
subspace induced from that of $\Bbb S^n_H$ gives rise to a
$k$-dimensional $\epsilon$-volume form and a $k$-dimensional
$\mu$-measure as a limit.

\begin{thm}\label{3.3}A $k$-dimensional region with piecewise analytic boundary contained in a $k$-dimensional
totally geodesic subspace transversal to $\partial\Bbb H^n$ has a
finite $k$-dimensional measure which is invariant under the
isometry action.
\end{thm}

\begin{proof}Let $U$ be a $k$-dimensional region contained in a
geodesic sphere $S^k$. Here we assume $S^k$ is transversal to
$\partial\Bbb H^n$ and we will discuss the case when $S^k$ is
tangent to $\partial\Bbb H^n$ in the following remark.

If $U$ does not intersect $\partial\Bbb H^n$, then the theorem is
clear and we assume $U$ intersect $\partial\Bbb H^n$. It suffices
to show that $\mu_k(U)$ is finite and $\mu_k(U)=\mu_k(gU)$ for any
$g\in PO(n,1)$ when $U$ is contained in $l^k:=\Bbb S^k_H=\Bbb
S^n_H\cap\{x|x_{k+1}=\cdots =x_n=0\}$. In this case $\mu_k(U)$ is
given as follows as in Proposition 3.2 and hence finite.
$$\al
\mu_k(U)&=\lim_{\epsilon\to 0} \int_U
dV^{l^k}_{\epsilon}=\lim_{\epsilon\to 0} \int^b_a
\frac{r^{k-1}F(r)}{(d^2_{\epsilon}-r^2)^{\frac{k+1}{2}}} dr\\
&=\lim_{\epsilon\to 0} \int_{\gamma}
\frac{r^{k-1}F(r)}{(d^2_{\epsilon}-r^2)^{\frac{k+1}{2}}}
dr=\int_{\gamma} \frac{r^{k-1}F(r)}{(1-r^2)^{\frac{k+1}{2}}}
dr=\int_{\gamma}\int_{B_r} dV_0^{l^k}, \eal
$$
when $dV^{l^k}_{\epsilon}$ is the $k$-dimensional
$\epsilon$-volume form on $l^k:=\Bbb S^k_H$.

\begin{figure}[h]
\begin{center}
\includegraphics[width=0.5\textwidth]{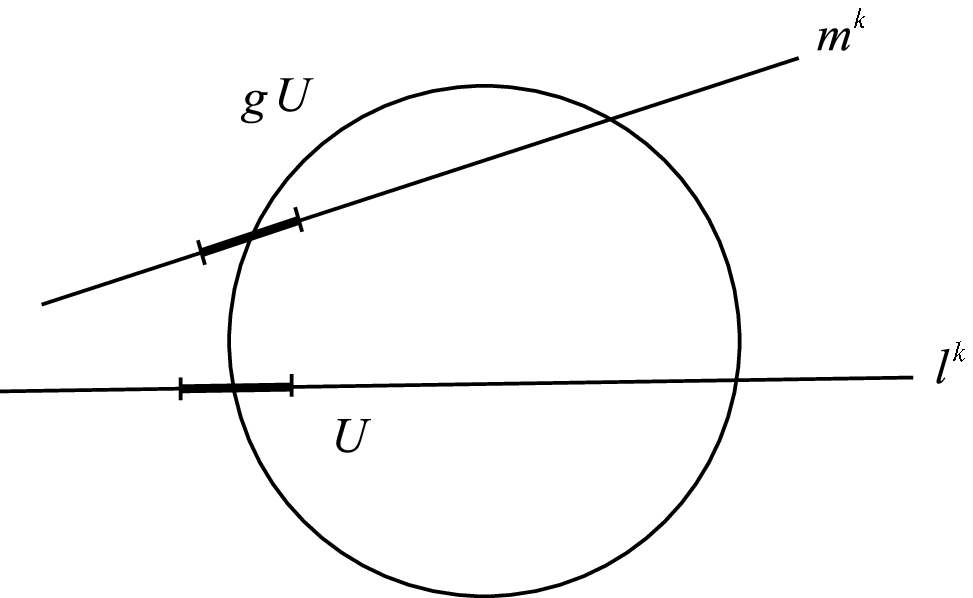}
\caption{\textbf{4.3}}
\end{center}
\end{figure}

Now we show that $\mu_k(gU)=\int_{\gamma}\int_{B_r} dV_0^{l^k}$
for $g\in PO(n,1)$ as before by checking the following steps.
$$\al
\mu_k(gU)&=\lim_{\epsilon\to 0} \int_{gU}
dV^{m^k}_{\epsilon}=\lim_{\epsilon\to 0} \int_{U}
g^*(dV^{m^k}_{\epsilon})=\lim_{\epsilon\to 0} \int^b_a \int_{B_r}
g^*(dV^{m^k}_{\epsilon})\\&=\lim_{\epsilon\to 0} \int_{\gamma}
\int_{B_r} g^*(dV^{m^k}_{\epsilon})=\int_{\gamma} \int_{B_r}
g^*(dV^{m^k}_0)=\int_{\gamma} \int_{B_r} dV_0^{l^k}.
 \eal
$$
But the proofs of these are exactly same if we notice that the
$k$-dimensional $\epsilon$-volume form $dV^{m^k}_{\epsilon}$ on
$m^k$ has the poles below the real $r$-axis near $r=1$. In fact
the poles of $dV^{m^k}_{\epsilon}$ are exactly the same as before,
i.e., $d^2_{\epsilon}-|x|^2=0$: Let
$e_1,\ldots,e_k,e_{k+1},\ldots,e_n$ be a local orthonormal frame
so that $\{e_1,\ldots,e_k\}$ spans the tangent space of $m^k$ and
$\alpha^1,\ldots,\alpha^k,\alpha^{k+1},\ldots,\alpha^n$ be its
dual frame so that $AE=I, A=(\alpha_{ij}), E=(e_{ij})$, where
$e_j=\sum e_{ij}\frac{\partial}{\partial x_i}$ and $\alpha^i=\sum
\alpha_{ij}dx^j$. If we denote the metric tensor $g$ by
$G=(g_{ij}), g_{ij}=g(\frac{\partial}{\partial
x_i},\frac{\partial}{\partial x_j})$, then $E^tGE=I$ and hence
$$\al
dV^n&=\alpha^1\wedge\cdots\wedge\alpha^n=\det A~~
dx^1\wedge\cdots\wedge dx^n=\frac{1}{\det E}~~
dx^1\wedge\cdots\wedge dx^n\\&=\sqrt{\det
G}~~dx^1\wedge\cdots\wedge
dx^n=\frac{d_{\epsilon}dx^1\wedge\cdots\wedge
dx^n}{(d^2_{\epsilon}-|x|^2)^{\frac{n+1}{2}}}. \eal
$$ Since $A=E^{-1}=\frac{1}{\det E}~\text{adj} E=\sqrt{\det G}~\text{adj} E$,
the only poles of $dV^k=\alpha^1\wedge\cdots\wedge\alpha^k$ are
those of $dV^n$ coming from $\sqrt{\det G}$.
\end{proof}

\begin{rem}\label{3.3a}
If we consider a geodesic sphere $S^{k}$ tangent to $\partial\Bbb
H^n$ and a region $U$ contained in $S^{k}$, we can show Theorem
\ref{3.3} holds for even more general type region $U$. In fact, if we
let $S^{k}$ be the $k$ dimensional subspace given by
$\{x_1,\ldots,x_k,0,\ldots,0,1\}$ which is tangent to
$\partial\Bbb H^n$ at $p=(0,\ldots,0,1)$ in $K^{n}$, then
$$dV_{\epsilon}^k=\frac{c_{\epsilon}~dx_1\wedge\cdots\wedge dx_k}
{(c^2_{\epsilon}-(x_1^2+\cdots+x_k^2))^{\frac{k+1}{2}}}=\frac{c_{\epsilon}r^{k-1}}
{(c^2_{\epsilon}-r^2)^{\frac{k+1}{2}}}~d\theta dr,
$$
where $c^2_{\epsilon}=d^2_{\epsilon}-1$. If $\epsilon\to 0$,
$dV_{\epsilon}^k$ becomes 0 outside the origin and the origin
itself becomes a pole. And it can be shown by direct computation
that $\mu_k(U)=0$ if $p\notin \overline U$ and $\mu_k(U)=\frac 12
\text{vol }(\Bbb S^k_H)$ if $p\in {\text Int }(U) $. Hence in this
case it is clearly invariant under isometry. This singular measure
$\mu_k$ has some interesting properties, which we do not want to
pursue in this paper, but we will state some of these which is related to
the theorem. The proofs are direct calculations. Even if $\mu_k$
is supported only at the origin, it is not Dirac measure since it
is only finitely additive. Indeed if $U$ is a region with
piecewise smooth boundary with origin at the boundary point, then
$\mu_k(U)$ is determined by the infinitesimal solid angle at the
origin since it can be shown easily that $\mu_k$ is dilation
invariant. This scale invariance also shows that $\mu_k$ is
invariant under similarity fixing the origin and hence $\mu_k$ is
invariant under isometry action since the derivatives of the
isometries fixing the origin are conformal at the origin.
\end{rem}

When we consider a $k$-dimensional totally geodesic region $U$, we
may assume that $U$ lies on $\Bbb S^{k+1}_H=\Bbb
S^n_H\cap\{x|x_{k+2}=x_{k+3}=\cdots =x_n=0\}$ by an isometry. Now
we have 3 cases as in the 1-dimensional case already discussed:
The $k$-dimensional subspace $S^k$ containing $U$ (i) meets
$\partial\Bbb H^{k+1}$ transversely, (ii) is tangent to
$\partial\Bbb H^{k+1}$, or (iii) does not meet $\partial\Bbb
H^{k+1}$.

For the case (i),  we may assume $U$ lies on $\Bbb S^{k}_H$, and
for the case (iii), $S^k$ is contained in $\Bbb L^{k+1}$ and we
may assume $U$ lies on the equator of $\Bbb S^{k+1}_H$ which is
essentially the same as the standard unit sphere $\Bbb S^k$ but
with $i^k$-factor for volume.

The case (ii), the $\mu_k$-measure of total space $S^k$ has value
$\text{vol }(\Bbb S^k_H)$ from Remark \ref{3.3a}. So we summarize
as follows.

\begin{pro}\label{3.4} The $k$-dimensional $\mu$-measure of any $k$-dimensional
 geodesic sphere in $\Bbb S^n_H$ is $\text{vol }(\Bbb S^k_H)$.
\end{pro}

\begin{rem}\label{3.5} For a curve $\gamma$ in $\Bbb S^n_H$
especially when it passes through $\partial\Bbb H^n$
transversally, it would be natural to define its length by
$$ \mu(\gamma)=\lim_{\epsilon \to 0}\int_{\gamma} ds_{\epsilon}.$$
Of course, we have to show that this definition is invariant under
isometry. We do not pursue this issue further in this paper.
\end{rem}

\begin{rem}\label{3.6}It is well known that we can use cross ratio
to define the distance between the two points (or the length of
geodesic line segment) in $\Bbb H^2$, and also can extend this to
$\Bbb S^2_H$ (see \cite{Sc}). But in this case there arises a subtle
and confusing choice problem for multi-valued logarithm.

\begin{figure}[h]
\begin{center}
\includegraphics[width=0.44\textwidth]{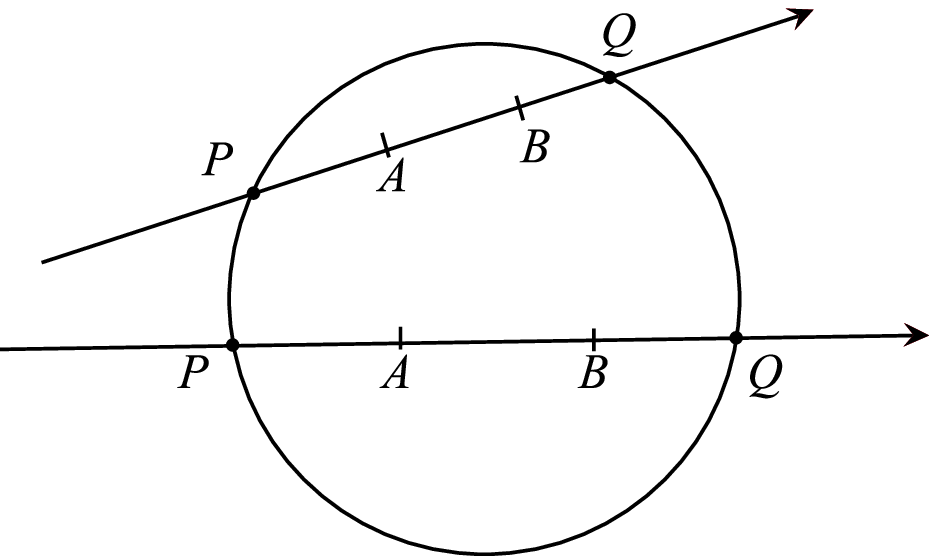}
\caption{\textbf{4.4}}
\end{center}
\end{figure}

For given two points $A$ and $B$ in $\Bbb H^2\subset K^2$, let $P$
and $Q$ be the points of intersection of $\partial\Bbb H^2$ and
the line passing through $A$ and $B$. We assign $P$ and $Q$ so
that $\overrightarrow{AB}$ and $\overrightarrow{PQ}$ have the same
direction. In this case the hyperbolic distance between $A$ and
$B$ is known to be as the following.
\begin{equation}\label{cro1}
d_H(A,B)=-\frac 12 \log(PQ|AB)=-\frac 12
\log\frac{\overline{PA}/\overline{AQ}}{\overline{PB}/\overline{BQ}}
\end{equation}
Here $\overline{PA},\cdots$ etc are Euclidean oriented length so
that $\overline{PA}=-\overline{AP}$ and $\overline{xy}$ is
positive when $\overrightarrow{xy}$ has the same direction as
$\overrightarrow{AB}$.

Now let's use the formula (\ref{cro1}) as a definition of distance
for any two points in $K^2$. Then in general one of the
multi-values of log in (\ref{cro1}) coincide with $\mu$-distance,
but it is unclear which choices are natural and consistent for all
cases.

Let $A=(x,0)$ and $B=(\frac 1x,0)$. Then $d_H(A,B)=-\frac 12
\log(-1)=\pm\frac {\pi}2i+m\pi i, ~~m\in\Bbb Z$. Taking the
clockwise contour from $1$ to $-1$, gives the value $\frac
{\pi}2i$, which coincides with $\mu$-distance of $A$ and $B$.

If the line passing through $A$ and $B$ is tangent to
$\partial\Bbb H^2$ at $P=Q$, we have to carry out a formal
computation for (\ref{cro1}) and choose log as follows to obtain
$\mu$-distance; $\log\frac{0/0}{1/-1}=\log(-1)=-\pi i$ or $\log
1=-2\pi i$ depending on whether $A=P=Q\ne B$ or $A\ne P=Q\ne B
~~(A<P<B)$, respectively.

The remaining case is when the line extending $AB$ does not meet
$\partial\Bbb H^2$. Let the equation of this line be $x_2=k
~~(k>1)$ and $A=(a,k), B=(b,k)$ with $a<b$ on $K^2$. If we solve
simultaneously with $x_1^2+x_2^2=1$, we have
$x_1=\pm\sqrt{k^2-1}i$. Now from (\ref{cro1}) with
$P=(\sqrt{k^2-1}i,k)$ and $Q=(-\sqrt{k^2-1}i,k)$ or
$P=(-\sqrt{k^2-1}i,k)$ and $Q=(\sqrt{k^2-1}i,k)$, we get
$d_H(A,B)=\pm\alpha i+m\pi i$ where $m\in\Bbb Z,~~0<\alpha<\pi$.
Among these, the choice $\alpha i$ coincides with the
$\mu$-distance and we should have chosen $P=(\sqrt{k^2-1}i,k)$ and
$Q=(-\sqrt{k^2-1}i,k)$ to obtain $\alpha i$, and it looks
confusing to explain why this is a consistent choice. The other
choice of $P$ and $Q$ corresponds to the choice $\tilde
d_{\epsilon}=1+\epsilon i$ for $\epsilon$-approximation of our
singular metric.
\end{rem}

The extended Kleinian model $\Bbb H^n_K$ has a projective
geometric structure, so a geodesic in $\Bbb H^n_K$ is a straight
line and a dual of a point $x$, i.e., $x^{\bot}$  is easily
obtained as usual (see Fig 4.5 and 4.6). Then the length of a
geodesic line segment joining $x$ (respectively $y$) and an
arbitrary point  in $x^{\bot}$ (respectively $y^{\bot}$) is
$\frac{\pi}2 i$. Because there is an isometry which takes  $x$ and
$x^{\bot}$ to a point on the equator and to a longitude
respectively, and takes $y$ and $y^{\bot}$ to a north pole and to
the equator respectively.

\begin{figure}[h]
\begin{center}
\includegraphics[width=0.3\textwidth]{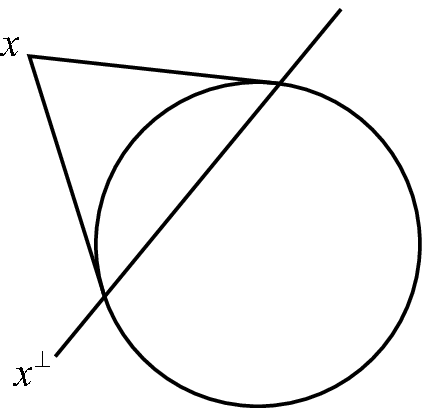}\quad\quad\quad
\includegraphics[width=0.35\textwidth]{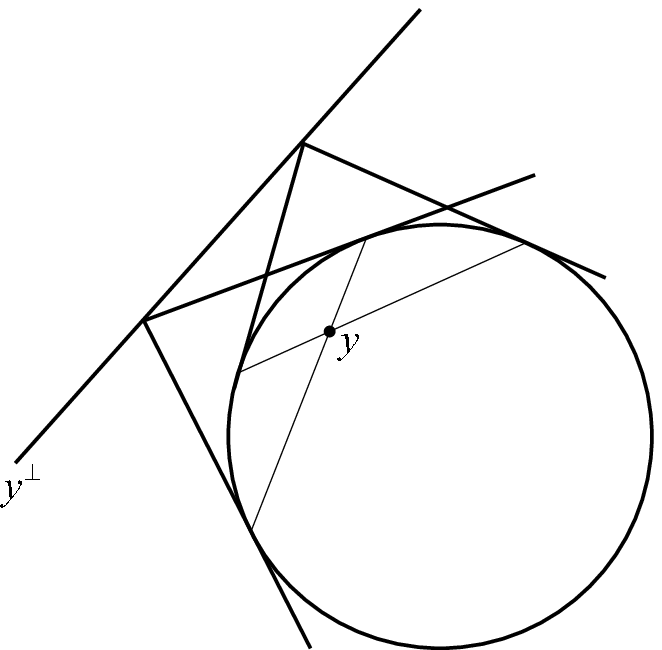}
\caption{\textbf{4.5}\textbf{\hskip 5truecm Fig 4.6\quad}}
\end{center}
\end{figure}

 Now we define angles on the extended model $\Bbb S_H^n$. In fact
we can define a notion of angle on a semi-Riemannian manifold
using $\mu$-measure on $\Bbb S_H^1$. For two tangent vectors $v_p$
and $w_p$ at a point $p$ on a Riemannian manifold, the notion of
the angle $\theta$ between these two vectors is obvious from that
of standard Euclidean plane $\Bbb R^2$ and can be calculated by
the equation,

\begin{equation}\label{6} \la v_p,w_p\ra= \vv v_p
\vv \vv w_p \vv \cos \theta,  \phantom{abcdef}  0 \leq \theta <
\pi.\end{equation}

But for a semi-Riemannian manifold, we have some difficulties with
this formula since the function $\cos^{-1}$ is multi-valued and
$\theta$ can take several complex values. The definitions of angle
have been given through the combinatorial way in \cite{D} and the
cross ratio in \cite{Sc}. Now that we have a notion of arc length on
$\Bbb S_H^1$, as $\mu$-measure, we can define an angle just as for
Riemannian case.

\begin{defi}\label{3.7}
For given two vectors $v,w \in \Bbb R^{n,1}$, the angle between
$v$ and $w$, $\theta = \angle(v,w)$, is defined as $-i\cdot
d_H(v,w)$, where $d_H(v,w)$ is the length of a geodesic segment
joining two points of $\Bbb S_H^n$  radially projected from $v,w$
to $\Bbb S_H^n$.
\end{defi}

Fig. 4.7 shows the various angles of between vectors in a
degenerate space which occurs as a tangent space to the light
cone. (Two poles in the picture represent light-like vectors of
$\Bbb R^{n,1}$.)

\begin{figure}[h]
\begin{center}
\includegraphics[width=0.3\textwidth]{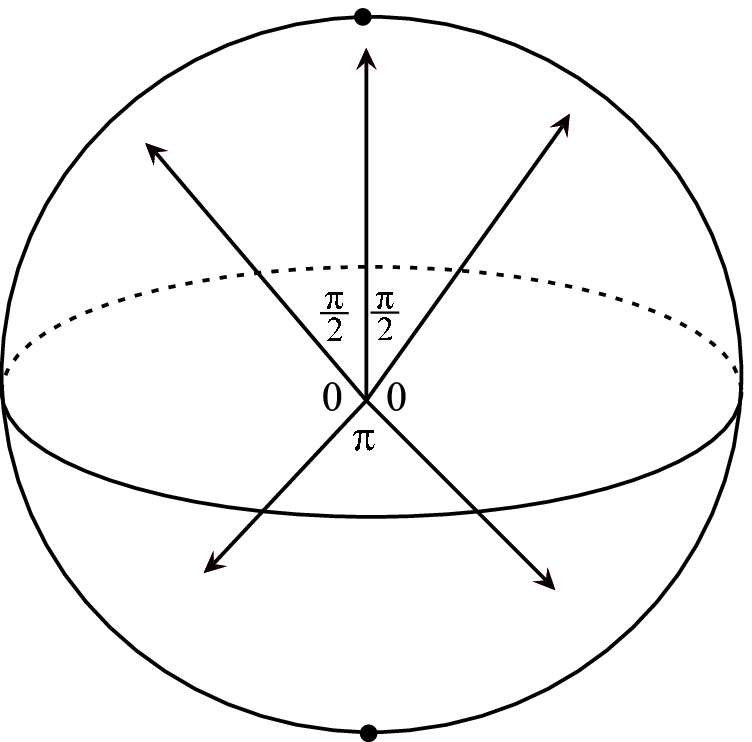}
\caption{\textbf{4.7}}
\end{center}
\end{figure}

Notice that the factor $-i$ is multiplied to normalize the total
length $2\pi i$ of the great circle in $\Bbb S_H^n$ as $2 \pi$
(see Theorem \ref{3.2}). If we used $\tilde d_{\epsilon} = 1 +
\epsilon i$ for our approximation, we have to use $i~d_H(v,w)$
instead of $-i~d_H(v,w)$ since the length of great circle in this
case becomes $-2\pi i$.

To check the equation (\ref{6}) for Lorentzian case, we need to
define the norm $\vv\cdot\vv$ first. The Lorentzian norm of a
vector $v$ in $\Bbb R^{n,1}$ is defined to be a complex number
$$ \vv v \vv = \la v,v\ra^{\frac{1}{2}},$$
where $\vv v\vv$ is either positive, zero, or positive imaginary.

From our definitions of angle and norm, it can be shown that
(\ref{6}) holds for all cases.

\begin{pro}\label{3.8}
For non-null vectors $v$ and $w$ in the space $\Bbb R^{n,1}$, we
have
\begin{equation}\label{7} \la v,w\ra= \vv v
\vv \vv w \vv \cos \angle(v,w).  \end{equation}
\end{pro}

\begin{proof}
It is not hard to check (\ref{7}) for the various cases through
simple computations, we can see the details at \cite{C}.
\end{proof}

Our definition of an angle clearly satisfies the following four
properties by Theorem \ref{3.3} and Remark \ref{3.3a},

\noindent (i) the invariance under isometry,\\
\noindent (ii) equation (\ref{7}),\\
\noindent (iii) finite additivity of angle: if $\theta$ consists
of two parts $\theta_1$ and $\theta_2$, then
$\theta=\theta_1+\theta_2$,\\
\noindent (iv) the angle of half rotation is $\pi$, i.e., a
straight line has angle $\pi$.

Conversely, it can be shown that the angle is uniquely determined
by these four properties.

There are various formulas similar to (\ref{7}) relating Lorentzian
inner product and angle or hyperbolic distance between two points or
its dual hyperplane depending on the position of vectors in
\cite[Proposition 2.4.5]{T},  and all of these are equivalent to
(\ref{7}) one single formula.

More generally if we consider the angle between two vectors $v$
and $w$ on $\Bbb R^{p,q}$ or two tangent vectors on $\Bbb S^n_H$
or on a semi-Riemannian manifold, we have to look at the plane
$P=span \{v,w\}$ spanned by these two vectors.

\begin{figure}[h]
\begin{center}
\includegraphics[width=0.35\textwidth]{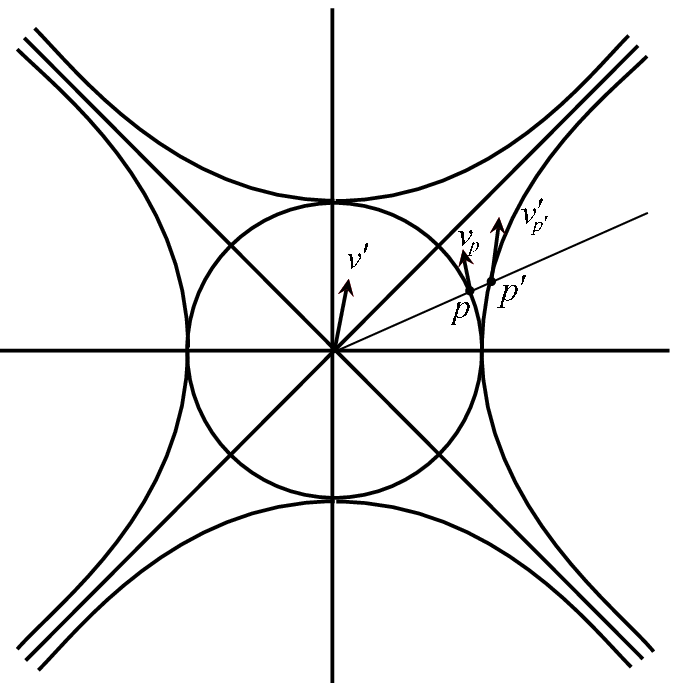}
\caption{\textbf{4.8}}
\end{center}
\end{figure}

The plane $P$ is isometric to $\Bbb {R}^{2,0}$ , $\Bbb {R}^{1,1}$,
or $\Bbb{R}^{0,2}$ for non-degenerate cases and to $\Bbb
{R}^{1,0}$, $\Bbb{R}^{0,1}$, or $\Bbb{R}^{0,0}$ for degenerate
cases depending on the situations. At a point $p$ on $\Bbb
S^n_{H}$ not lying on $\partial \Bbb {H}^n$, the angle between the
two tangent vectors $v_p$ and $w_p$ is determined on the plane
$span \{ v'_{p'}, w'_{p'} \}$ where $p'$ is the radial projection
of $p$ on the hyperboloid $H^n_{\pm}$ or $S^n_1$ and $v'_{p'}$,
$w'_{p'}$
   $\in T_{p'} ( H^n_\pm ) $
   or $T_{p'} (S^n_1 )$ are projections of $v_p$,
$w_p$ respectively. The plane $span \{ v'_{p'}$ , $w'_{p'}  \}$
has signature $\Bbb {R}^{2,0}$,$\Bbb {R}^{1,1}$  or $\Bbb
{R}^{1,0}$ depending on the position of $p$, and in each case the
angle is determined in the usual way for the definite cases and
using Definition \ref{3.7} for $\Bbb R^{1,1}$ and $\Bbb R^{1,0}$
cases.

If $p\in\Bbb S^n_H$ lies in the Lorentzian part $S^n_1$, the norm
$\vv v_p\vv$ of a tangent vector $v_p$ and the norm of its
translate to the origin, denoted by $v'$, are related by, $\vv
v'\vv=-i\vv v_p\vv$ from our sign conventions (see Fig. 4.8). And
$\la v',w'\ra=-\la v_p,w_p\ra$ by the definition of the metric on
$S^n_1$. Hence we have the following identity
\begin{equation}\label{x}
\frac{\la v_p,w_p\ra}{\vv v_p\vv \vv w_p\vv}=\frac{\la v',w'\ra}{\vv
v'\vv \vv w'\vv},\end{equation} and $\angle(v_p,w_p)=\angle(v',w')$
by the four characterizing properties of angle.

\begin{pro}\label{3.9} For $p\in \Bbb S^n_H$ and $p\notin\partial\Bbb
H^n$, the angle $\angle(v_p,w_p)$ between two tangent vectors
$v_p,w_p\in T_p\Bbb S^n_H$ is equal to $-i\dot d_H(v',w')$, where
$v'$ is a vector at the origin of $\Bbb R^{n,1}$ given as a
parallel translation of the tangent vector of $T_p H^n_{\pm}$ or
$T_p S^n_1$ which is obtained as a radial projection of $v_p$. (We
also can see the point of $K^n$ determined by $v'$ as in Fig.
4.9.)
\end{pro}

\begin{proof}If $p$ is in the Lorentzian part, we may assume $p$
lies on the equator via isometry. In this case, $v'=v$ and $w'=w$,
and  $\angle(v_p,w_p)=\angle(v,w)=-i~d_H(v,w)=-i~d_H(v',w')$ hold
from (\ref{x}) and Definition \ref{3.7}.

If $p$ is in the hyperbolic part, we may assume $p$ is the north
pole and we clearly have $d_H(v',w')=d_H(v,w)=i\angle(v_p,w_p)$.
\end{proof}

The isometry invariance of angle at a point $p\in \Bbb
S^n_H\backslash \partial\Bbb H^n$ is obtained from Proposition
\ref{3.9}.

\begin{figure}[h]
\begin{center}
\includegraphics[width=0.5\textwidth]{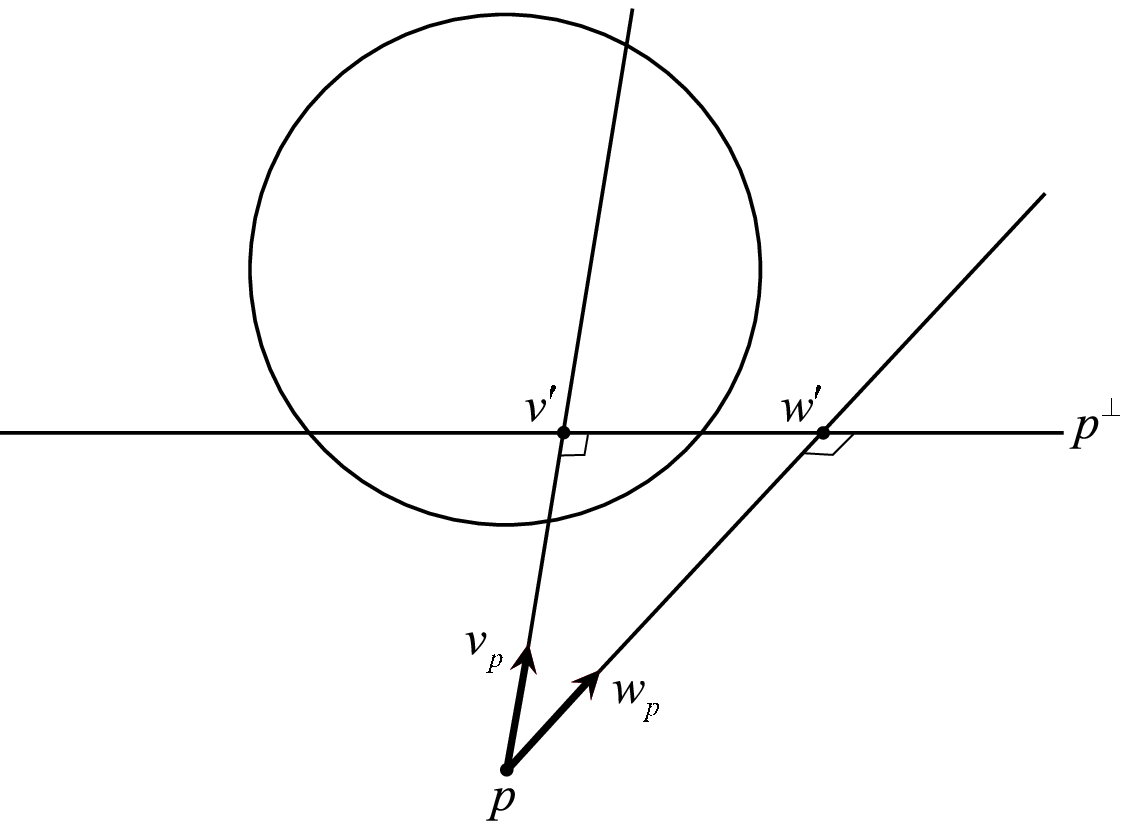}
\caption{\textbf{4.9}}
\end{center}
\end{figure}

If $p\in\Bbb S^2_H$ lies on $\partial\Bbb H^2$, we define angle as
follows.

\begin{defi}\label{3.10}We define the angle $\angle(v_i,v_{i+1})$
for the tangent vectors $v_1,\ldots,v_5$ configured as in Fig.
4.10:

\noindent $$\angle(v_1,v_2)=\angle(v_5,v_1)=\frac {\pi}{2},~~
\angle(v_2,v_3)=\angle(v_4,v_5)=0, \text{ and
}\angle(v_3,v_4)=\pi.$$
\end{defi}

\begin{figure}[h]
\begin{center}
\includegraphics[width=0.35\textwidth]{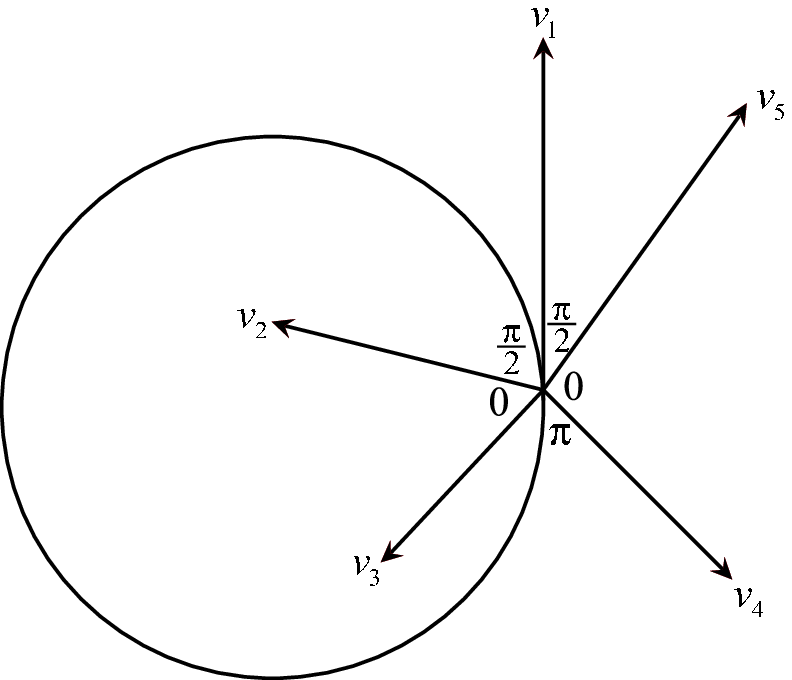}
\caption{\textbf{4.10}}
\end{center}
\end{figure}

\begin{rem}\label{3.11}A justification of this definition can be
given as follows. For $p\in\partial\Bbb H^2$  in $K^2$, we first
consider $\theta_{\epsilon}$ from the equation (\ref{6}) using
$\epsilon$-metric and then determine $\theta$ as $\lim_{\epsilon
\rightarrow 0} \theta_{\epsilon}$. In this process we again have a
problem of choice of values for multi-valued function $\cos^{-1}$.
But our choice is the simplest reasonable  one satisfying the four
characterizing of angle.

If $p \in \partial \Bbb H^n$ with $n \geq 3 $, we have another
type of tangent plane at $p$ which touches $\partial \Bbb H^n $ at
the only point $p$. In this case, the $\epsilon$-metric and hence
the limit has a rotational symmetry around $p$ and hence the angle
is the same as the usual Euclidean angle.
\end{rem}

\begin{rem}\label{3.12}When we consider an angle on $\Bbb S^n_H$,
we have 3 different approaches other than using (Definition
\ref{3.7}, \ref{3.10}) or (Proposition \ref{3.9}, Definition
\ref{3.10}) or (four characterizing properties of angle).

One way is to use the dual distance and we will defer the discussion
of this topic to elsewhere.

Secondly  we can define the angle between $v_p$ and $w_p$ as the
area of the lune determined by these two vectors. The angle is
given by $-\frac 12 \times$ the area (2-dimensional $\mu$-measure)
of the lune.

Third definition uses cross ratio and will be explained separately
in the next remark.
\end{rem}

\begin{rem}\label{3.13}As we explained distance on $K^2$ using
cross ratio, we also can explain angle using cross ratio. Let $P$ be
a point in the Lorentzian part, and $\theta$ be an angle between two
lines $a$ and $b$. Let $s$ and $t$ be two lines through $p$ tangent
to $\partial\Bbb H^2$, and be ordered so that the directions
$\overrightarrow{ab}$ and $\overrightarrow{st}$ match.

\begin{figure}[h]
\begin{center}
\includegraphics[width=0.45\textwidth]{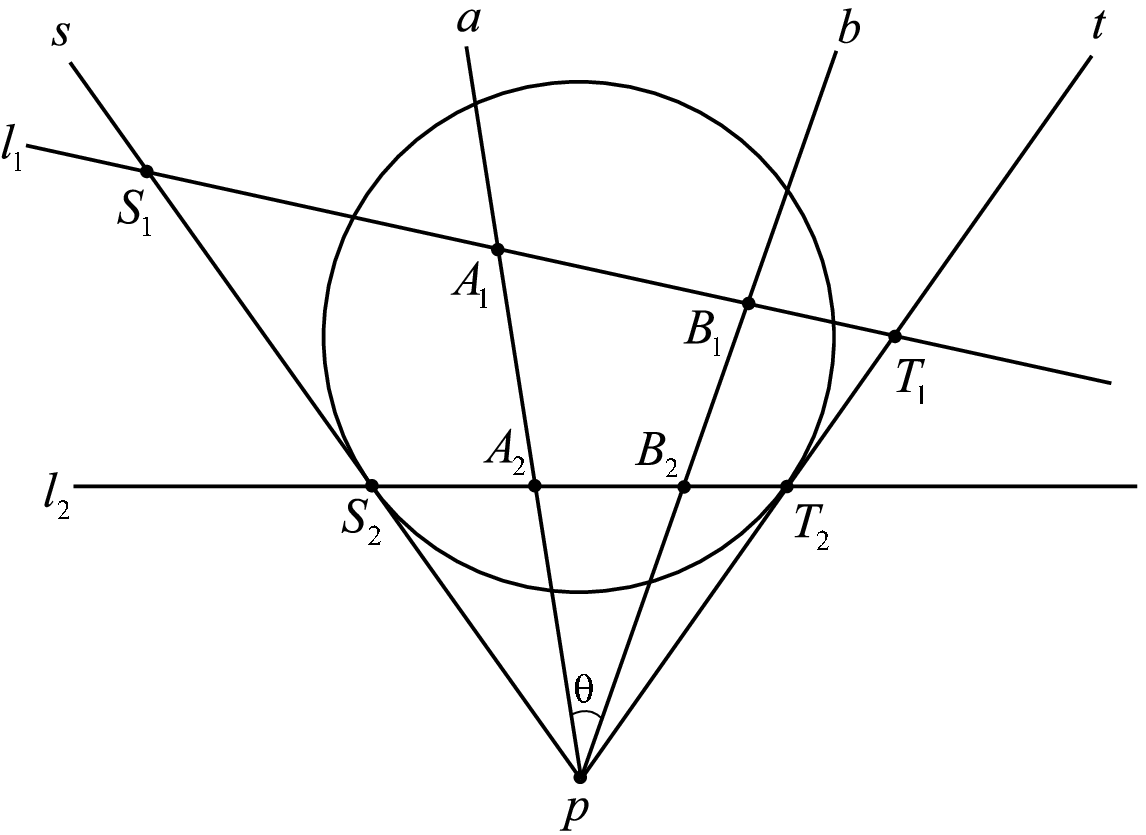}
\caption{\textbf{4.11}}
\end{center}
\end{figure}

Choose any two lines $l_i, i=l,2$ and let $S_i,A_i,B_i,T_i$ be the
intersecting points of $l_i$ and the lines $s,a,b,t$ respectively.
Then it is well known from a Pappus theorem that the cross ratio
$(S_i T_i|A_i B_i)$ agrees for $i=1,2$, and hence we can define
$(st|ab):=(S_i T_i|A_i B_i)$ and the angle $\theta$ is defined as
$\angle\theta=\frac{i}{2}\log (st|ab)$ (see \cite{Sc}).

Now it can be shown using Remark \ref{3.6} and Proposition
\ref{3.9} that this definition agrees with ours.

If $p$ lies inside $\Bbb H^2$, we can compute the tangent lines $s$
and $t$ through complexification. In this case, we have some
ambiguity for the choice of $s$ and $t$ as in the later half of
Remark \ref{3.6}. Anyway Remark \ref{3.6} and Proposition \ref{3.9}
also imply that our angle $\theta$ agrees with an appropriate choice
of angle comming from the cross ratio.

When $p$ lies on $\partial\Bbb H^2$, the two tangent lines $s$ and
$t$ coincide and we can derive the same value as in the Fig. 4.10
from the cross ratio by choosing an appropriate values of log.
\end{rem}

\section{Gauss volume formula and Santal\'o-Milnor-Su\'arez-Peir\'o relation
in the extended hyperbolic model}

Choose a point $x\in \Bbb S^n_H$ and consider its dual hypersphere
$x^{\bot}\cap \Bbb S^n_H$. For a subset $B\subset x^{\bot}\cap \Bbb
S^n_H$, the union of the great semi-circles passing through $b\in B$
connecting $x$ and $-x$ will be called a lune $L=L(x,B)$ with vertex
$x$ and cross section $B$. And the ($n-1$)-dimensional $\mu$-measure
$\mu(B)$ in $x^{\bot}\cap \Bbb S^n_H$ will be called a hyperbolic
solid angle of $L$. In fact, if $x=e_0=(1,0,\cdots)\in\Bbb R^{n,1}$,
then the dual hypersphere is $S^n_1\cap\{x\in\Bbb R^{n,1}|x_0=0\}$
and the radial projection (or an inverse image of exponential map)
of a lune $L(e_0,B)$ becomes a cone in the tangent space
$T_{e_0}(\Bbb S^n_H)$ and $\mu(B)$ is $i^{n-1}$ times the usual
solid angle of the cone in Euclidean space $T_{e_0}(\Bbb S^n_H)$,
i.e., $i^{n-1}$ times spherical volume of the portion in $\Bbb
S^{n-1}$ cut out by the cone. We should be careful when we consider
the volume of $\Bbb S^{n-1}$ in $\Bbb S^n_H$ since it lies in $\Bbb
L^n$ and its metric is the negative of the usual one and the volume
is determined as $i^{n-1}\cdot\text{vol }(\Bbb S^{n-1})$ by
Convention and of course is equal to $\text{vol }(\Bbb S^{n-1}_H)$.

\begin{thm}\label{4.1} A lune of hyperbolic solid angle $S$ in $\Bbb
S^n_H$ has volume $\frac{\text{vol }(\Bbb S^n_H)}{\text{vol }(\Bbb
S^{n-1}_H)}\cdot S $, where the vertex of the lune assumed not to
lie on $\partial\Bbb H^n$ for convenience. (See the following remark
for the case when the vertex lies on $\partial\Bbb H^n$.)
\end{thm}
\begin{proof} Let's consider first the case when the vertex of a
lune lies in $\Bbb H^n$. In this case we may assume the vertex is
$e_0=(1,0,\cdots)$, i.e., vertex is $(0,\cdots,0)\in K^n$,  via
isometry, and the volume of a lune is clearly proportional to the
solid angle by the rotational symmetry of $\Bbb S^n_H$.

When the vertex of a lune lies in the Lorentzian part $S^n_1$, we
may assume the vertex is $e_1=(0,1,0,\cdots)$, i.e., vertex is
$(\infty,0,\cdots,0)\in K^n$, by an isometry. Then such lune $L$
can be given as a cylinder of type $(-\infty,\infty)\times S$ in
the Kleinian model, and the volume of $L$ will be given by
$$\mu(L)=\lim_{\epsilon \to 0}\int_S \int_{-\infty}^{\infty}
\frac{d_{\epsilon}dx_1}{(d_{\epsilon}^2 -x_1^2-x_2^2-\cdots
-x_n^2)^{\frac{n+1}2}} ~~dx_2 \cdots dx_n \quad\quad (n\ge 2).$$
To compute the integral, we divide $S$ into two parts $A$ and $B$
depending on whether $x_1^2+x_2^2+\cdots +x_n^2$ is greater than 1
and less than or equal to 1 respectively. Since the sign of the
volume of Lorentzian part is $i^{n+1}$, the integral
$\int_{-\infty}^{\infty} \frac{d_{\epsilon}dx_1}{(d_{\epsilon}^2
-x_1^2-x_2^2-\cdots -x_n^2)^{\frac{n+1}2}}$ in the above formula
for $A$ part can be replaced by $i^{n+1}\int_{-\infty}^{\infty}
\frac{d_{\epsilon}dx_1}{(x_1^2+x_2^2+\cdots
+x_n^2-d_{\epsilon}^2)^{\frac{n+1}2}}$ when the integrand has
positive real part (and positive as $\epsilon \to 0$). If we let
$a=x_2^2+\cdots +x_n^2-d_{\epsilon}^2$ and $x_1=a^{\frac 12}t$
with $Re(a^{\frac 12})>0$, then
$$\al
i^{n+1}\int_{-\infty}^{\infty}
\frac{dx_1}{(x_1^2+a)^{\frac{n+1}2}}&=\frac{i^{n+1}}{a^{\frac
n2}}\int_{\gamma_1} \frac{2~ dt} {(1+t^2)^{\frac{n+1}2}}=\frac{2
i^{n+1}}{a^{\frac n2}}\int_0^{\infty}
\frac{dt}{(1+t^2)^{\frac{n+1}2}}\\&=\frac{2i}{(d_{\epsilon}^2
-x_2^2-\cdots -x_n^2)^{\frac n2}}\int_{0}^{\infty}
\frac{dt}{(1+t^2)^{\frac{n+1}2}}.\eal$$ Here
$\frac{1}{(d_{\epsilon}^2 -x_2^2-\cdots -x_n^2)^{\frac
n2}}=\frac{i^{n}}{a^{\frac n2}}$. Since $A$ is a Lorentzian part
in $K^{n-1}$ and $\gamma_1$ is a contour given in Fig. 5.1.

For the $B$ part, let $b=d_{\epsilon}^2-x_2^2-\cdots -x_n^2$ and
$x_1=b^{\frac 12}y$ with $Re(b^{\frac 12})>0$, then
$$
\int_{-\infty}^{\infty}
\frac{dx_1}{(b-x_1^2)^{\frac{n+1}2}}=\frac{2}{b^{\frac
n2}}\int_{\gamma_2}
\frac{dy}{(1-y^2)^{\frac{n+1}2}}=\frac{2}{b^{\frac
n2}}\int_{\gamma_3} \frac{dy}{(1-y^2)^{\frac{n+1}2}}.$$
Substituting $y=i t$, we obtain
$$
\frac{2}{b^{\frac n2}}\int_{\gamma_3}
\frac{dy}{(1-y^2)^{\frac{n+1}2}}=\frac{2i}{b^{\frac
n2}}\int_{\gamma_4}
\frac{dt}{(1+t^2)^{\frac{n+1}2}}=\frac{2i}{b^{\frac
n2}}\int_{0}^{\infty} \frac{dt}{(1+t^2)^{\frac{n+1}2}}.$$

\begin{figure}[h]
\begin{center}
\includegraphics[width=0.9\textwidth]{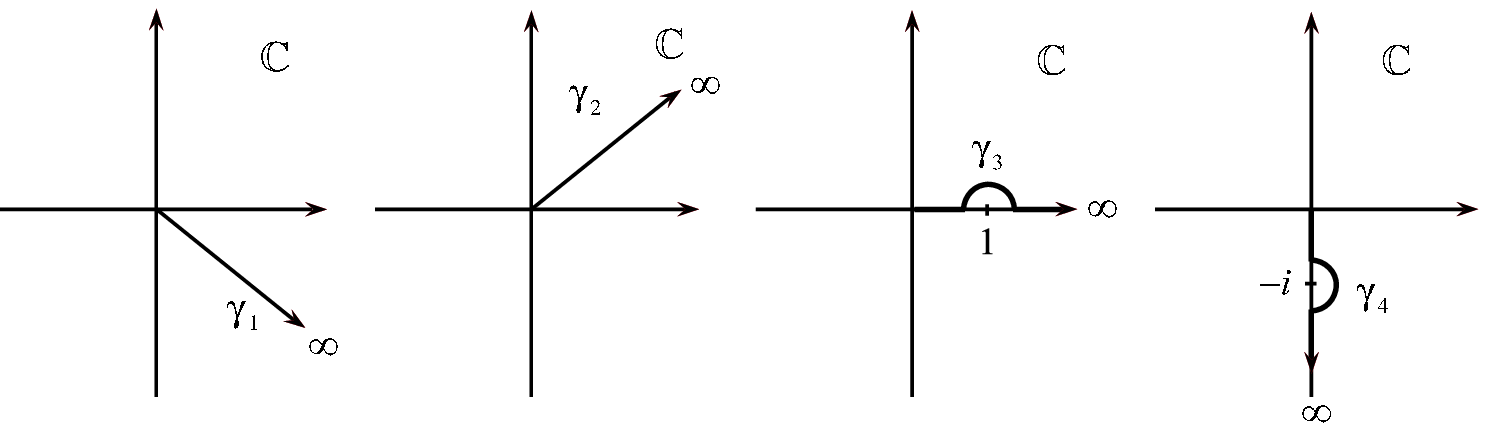}
\caption{\textbf{5.1}}
\end{center}
\end{figure}

{\noindent\bf Claim.} $\int_{0}^{\infty}
\frac{dt}{(1+t^2)^{\frac{n+1}2}}=\frac{\text{vol }(\Bbb
S^n)}{2\text{vol }(\Bbb S^{n-1})} \quad\quad (n\ge 2).$

Now if assume this claim, we can write for both $A$ and $B$ parts
as
$$
\mu_n(L)=\lim_{\epsilon \to 0}\frac{i\text{vol }(\Bbb
S^n)}{\text{vol }(\Bbb S^{n-1})}\int_{S}
\frac{d_{\epsilon}dx_2\cdots dx_n}{(d_{\epsilon}^2 -x_2^2-\cdots
-x_n^2)^{\frac{n}2}}=\frac{\text{vol }(\Bbb S^n_H)}{\text{vol
}(\Bbb S^{n-1}_H)}\mu_{n-1}(S),
$$
where $\mu_n(L)$ is the $n$-dimensional volume or $\mu$-measure.
And this completes the proof of the theorem.

{\noindent\bf Proof of Claim} Let $I_m=\int_{0}^{\infty}
\frac{dt}{(1+t^2)^m}$. Then $I_m=\frac{2m-3}{2m-2}I_{m-1}~(m\ge
2)$ is induced from
$$\int\frac{dt}{(1+t^2)^m}=\frac{1}{2(m-1)}\frac{t}{(1+t^2)^{m-1}}+\frac{2m-3}{2m-2}\int\frac{dt}{(1+t^2)^{m-1}}.$$
Prove the Claim using induction on $n\ge 2$.
 {\noindent If $n=2,
I_{3/2}=[\frac{t}{\sqrt{1+t^2}}]_0^{\infty}=1=\frac{\text{vol
}(\Bbb S^2)}{2\text{vol }(\Bbb S^1)}.$}

{\noindent And if $n=3$, $I_2=\frac 12 I_1=\frac
{\pi}{4}=\frac{\text{vol }(\Bbb S^3)}{2\text{vol }(\Bbb S^2)}.$}
Now using the hypothesis, we have
$$I_{\frac{n+2}{2}}=\frac{n-1}{n}I_{\frac n2}=\frac{n-1}{n}\frac{\text{vol
}(\Bbb S^{n-1})}{2\text{vol }(\Bbb S^{n-2})}=\frac{\text{vol
}(\Bbb S^{n+1})}{2\text{vol }(\Bbb S^n)}.$$
The last equality
follows from the well-known formula
$$\frac{\text{vol
}(\Bbb S^{n})}{\text{vol }(\Bbb S^{n-2})}=\frac{2\pi}{n-1}\quad
\text{or} \quad \text{vol }(\Bbb
S^n)=\frac{2\pi^{\frac{n+1}2}}{\Gamma(\frac{n+1}2)}.
$$
\end{proof}

\begin{rem}\label{r1} The above theorem seems still hold even when
the vertex of a lune lies on the ideal boundary if the lune is
transversal to the ideal boundary. The non-transversal case looks
rather complicated and subtle. We illustrate these when the
dimension $n=2$ to show some ideas.

\begin{figure}[h]
\begin{center}
\includegraphics[width=0.7\textwidth]{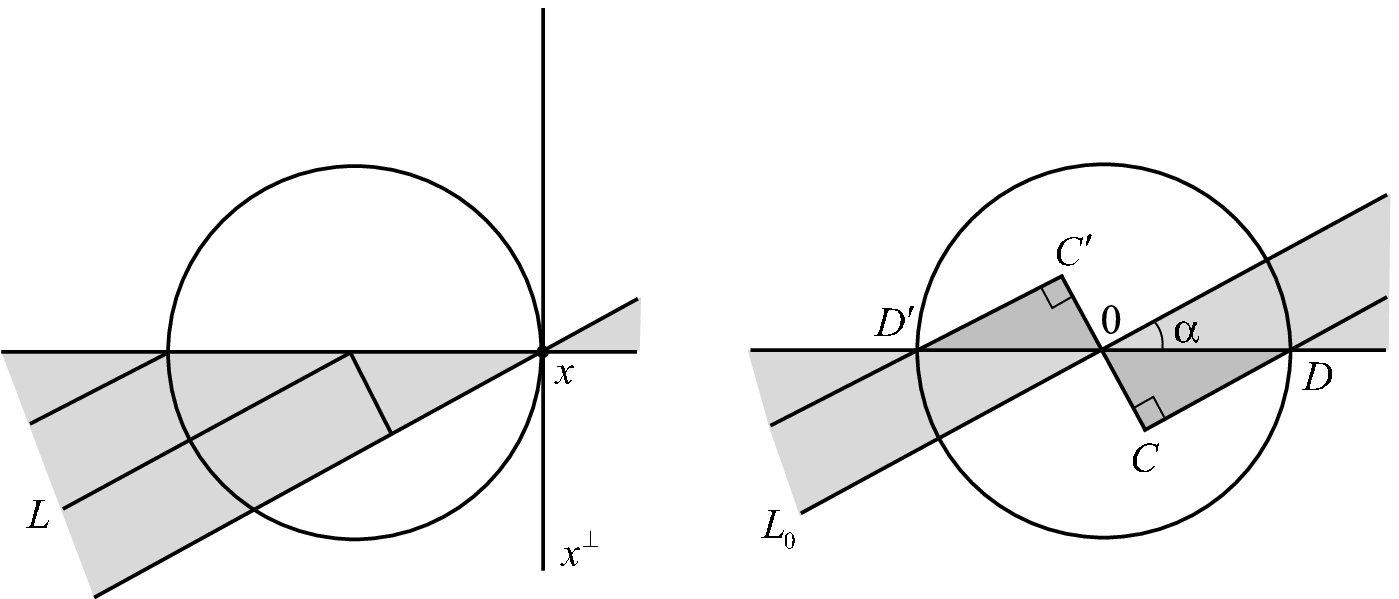}
\caption{\textbf{5.2}}
\end{center}
\end{figure}

The first picture of Figure 5.2 is a lune $L$ with its vertex
$x\in \partial\Bbb H^2$ and its angle $0$. We can change $L$ into
a region consisting of a lune $L_0$ with its vertex at the origin
and angle $\alpha$ and two right triangles as in the second
picture so that it clearly has the same area as $L$. Now from
this, we can show the area of $L$ is in fact zero as follows.
$$
\al
\mu(L)&=\lim_{\epsilon \to 0}\int_L dV_{\epsilon}\\
 &=\lim_{\epsilon \to 0}\int_{\triangle OCD\cup\triangle OC'D'\cup L_0} dV_{\epsilon} \\
 &=2\mu(\triangle OCD)+\mu(L_0)\\
 &=2(\pi-\frac{\pi}2-(\frac{\pi}2-\alpha)-0)+(\frac{-4\pi}{2\pi i})\alpha i\\
 &=0.
\eal
$$

On the other hand, for this $L$, $x^{\bot}\cap L=\{x\}$ and hence
$\mu(B)=0$. Therefore we see that the theorem holds for such $L$.

\begin{figure}[h]
\begin{center}
\includegraphics[width=0.7\textwidth]{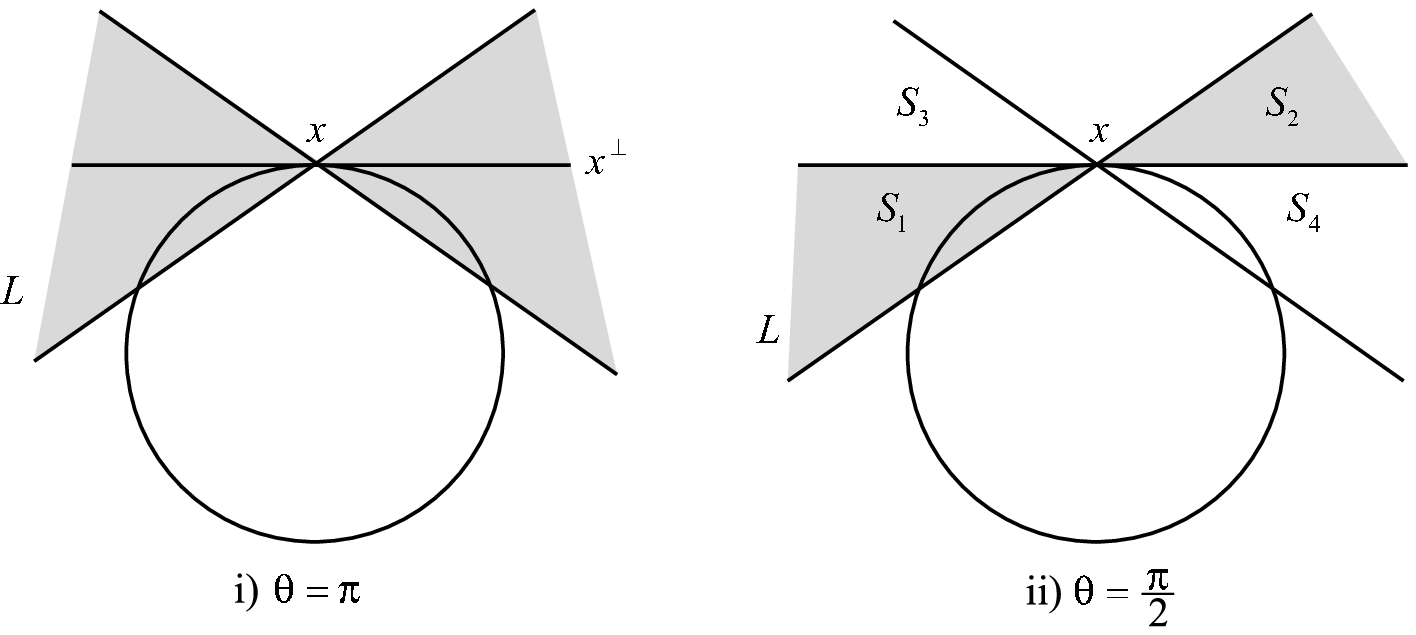}
\caption{\textbf{5.3}}
\end{center}
\end{figure}

If a Lune $L$ has an angle $\pi$ as in the figure 5.3 (i), then
its complement $K^2\backslash L$ is a lune of angle $0$. Hence we
have
$$
\mu(L)=-2\pi=\frac{-4\pi}{2\pi i}\cdot \pi i=\frac{\text{vol
}(\Bbb S^2_H)}{\text{vol }(\Bbb S^1_H)} \mu_1(x^{\bot}\cap L)
$$
and the theorem holds.

If a lune $L$ has an angle $\frac{\pi}2$ as in figure 5.3 (ii),
then let $S_3$ and $S_4$ are the reflection of $S_2$ and $S_1$
respectively about the vertical line through $x$. Even if
$\mu(S_1)$ and $\mu(S_2)$ are not defined separately, we can
compute $\mu(L)$ as follows using the existence of $\mu(S_1\cup
S_3)$.
$$
\al \mu(L)&=\lim_{\epsilon \to 0}\int_{S_1\cup S_2} dV_{\epsilon}=\lim_{\epsilon \to 0}(\int_{S_1} dV_{\epsilon}+
\int_{S_2}
dV_{\epsilon}) \\
 &=\lim_{\epsilon \to 0}(\int_{S_1} dV_{\epsilon}+\int_{S_3} dV_{\epsilon})=\lim_{\epsilon \to 0}\int_{S_1\cup S_3}
 dV_{\epsilon}\\
 &=\mu(S_1\cup
S_3)=(\mu(S_1\cup S_3)+\mu(S_2\cup S_4))/2\\
 &=-\pi.
\eal
$$
In this case, it is unclear how to define a solid angle. But the
above result suggests that it would be reasonable to interpret
that $L$ occupies the half of the total solid angle.
\end{rem}

\begin{rem}\label{r2} When the dimension $n=2$, we have to be
careful about the notions of angle and solid angle. If the angle
of a lune is $\theta$, then its solid angle becomes $\theta i$
from its definition.
\end{rem}

When the dimension $n$ is even, the volume of a spherical simplex on
$\Bbb S^n$ can be obtained combinatorially as an alternating sum of
its solid angle using the well known Euler-Poincar\'e method (see
\cite[p.120]{Geo}). And then it can be shown that the same formula
(only differ by sign) holds for  a hyperbolic simplex indirectly
using an analytic continuation technique. But if we use the extended
model $\Bbb S^n_H$, the above Euler-Poincar\'e method can be applied
directly without any change showing that the above formula holds not
only for a hyperbolic simplex but also for a Lorentzian simplex or
even for a simplex lying across the ideal boundary $\partial \Bbb
H^n$, that is, any simplex on $\Bbb S^n_H$.

When we consider the formula, it is more convenient to use the
normalized volume so that the  total volume becomes 1. The
normalized volume of a lune $L$ is given by
$\widehat{\mu}(L)=\frac{\mu(L)}{\text{vol }(\Bbb S^n_H)}$ and the
normalized solid angle of $L$ is given by
$\widehat{\alpha}(L)=\frac{\alpha(L)}{\text{vol }(\Bbb
S^{n-1}_H)}$.

Given an $n$-simplex $\triangle^n$ on $\Bbb S^n_H$, $\triangle^n$
is an intersection of  half spaces $H_0, H_1, \ldots , H_n$ and
$H_{i_1}\cap\cdots\cap H_{i_k}$ $(k\leq n)$ becomes a lune
$L_{i_1\ldots i_k}$. Then the normalized volume of $\triangle^n$
can be given as follows using Euler-Poincar\'e method and Theorem
5.1.

\begin{equation}\label{B}
\al
\widehat{\mu}(\triangle^n)&=\frac12\sum_{k=0}^n\sum_{i_1<\cdots<i_k} (-1)^k \widehat{\mu}(L_{i_1\cdots i_k})\\
&=\frac12  \sum_{k=0}^n\sum_{i_1<\cdots<i_k} (-1)^k \widehat{\alpha}(L_{i_1\cdots i_k})\\
&=:\frac12 k(\triangle^n).
\eal
\end{equation}

If we apply the formula (\ref{B}) when $n=2$ for a triangle
$\triangle$ on $\Bbb S^2_H$ with three angles $A, B,$ and $C$,
then we have
$$
\frac{\mu(\triangle)}{-4\pi}=\widehat{\mu}(\triangle)=\frac12(1-\frac12-\frac12-\frac12+\frac{Ai}{2\pi
i}+\frac{Bi}{2\pi i}+\frac{Ci}{2\pi i}).
$$

It follows that $\mu(\triangle)=\pi-A-B-C$, and we see that this
beautiful formula can be extended even across the ideal boundary
$\partial\Bbb H^2$. This was first obtained by J. B\"ohm and Im
Hof \cite{I}.

Suppose $M$ is a hyperbolic manifold with a triangulation $\mathcal
T$ consisting of totally geodesic  $n$-simplex. Then it can be shown
easily that $\chi(M)=\sum_{\triangle^n\in\mathcal T}k(\triangle^n)$
(see \cite{JK}). Therefore the Gauss-Bonnet theorem for a hyperbolic
manifold follows immediately.

\begin{equation}\label{A}
\al
\mu(M^n)=\text{vol }(M^n)&=\text{vol }(\Bbb S^n_H)\sum_{\triangle^n\in\mathcal T}\widehat{\text{vol}}(\triangle^n)=
\text{vol }(\Bbb S^n_H)\sum_{\triangle^n\in\mathcal T}\widehat{\mu}(\triangle^n)\\
&=\frac{\text{vol }(\Bbb S^n_H)}2\sum_{\triangle^n\in\mathcal T} k(\triangle^n)=\frac{\text{vol }(\Bbb S^n_H)}2\chi(M^n).
\eal
\end{equation}

Therefore we conclude the following result from the
Euler-Poincar\'e method used in $\Bbb S^n_H$.
\begin{pro}\label{4.2}When the dimension $n$ is even, we have\\ $\widehat{\mu}(\triangle^n)=\frac12
\sum_{k=0}^n\sum_{i_1<\cdots<i_k} (-1)^k \widehat{\alpha}
(L_{i_1\cdots i_k})$ for an $n$-simplex $\triangle^n$ in $\Bbb
S^n_H$ and \\$\mu(M^n)=\frac{\text{vol }(\Bbb S^n_H)}2\chi(M^n)$
for a  hyperbolic $n$-manifold. In particular, the area of a
triangle in $\Bbb S^2_H$ is $\pi-A-B-C$, where  $A,B,C$ are the
angles of the given triangle.
\end{pro}

Similarly if we consider a Lorentzian spherical manifold $M$
(i.e., metric signature $(-,$ $+,$ $\cdots,$ $+)$ with constant
sectional curvature $\mathcal K\equiv 1$, or metric $(+,$ $-,$
$\cdots,-)$ with $\mathcal K\equiv -1$), then $M$ has a developing
on the Lorentzian part and we have the same formula
$\mu(M^n)=\frac{\text{vol }(\Bbb S^n_H)}2\chi(M^n)$.

In this argument we have to be careful not to have a simplex
$\triangle^n$ whose (extended) face is tangent to $\partial\Bbb H^n$
since the solid angles are not defined for such a simplex. But of
course this can be easily achieved by perturbing the triangulation.

Now notice that $\mu(M^n)$ has sign $-i^{n-1}=i^{n+1}$ and
$\text{vol }(\Bbb S^n_H)$ has sign $i^n$. This shows that we have
both $\mu(M^n)=0$ and $\chi(M^n)=0$. We already know that latter
should hold since $M$ is Lorentzian, but the condition $\mu(M^n)=0$
is absurd and we can conclude that there does not exist such closed
manifold $M$. Of course, this fact has been known and can be deduced
from the usual Gauss-Bonnet theorem for semi-Riemannian manifolds
(see \cite{Ku} for instance), but we could see this immediately by
an elementary combinatorial way.

Proposition \ref{4.2} gives us an interesting consequences
especially when the $\triangle(A,B,C)$ lies across the ideal
boundary $\partial\Bbb H^2$ such that the edges are not completely
contained in the Lorentzian part. In this case $\pi-A-B-C$ is a
complex number whose real part is the area of truncated polygon
and whose imaginary part is the length of the edge introduced in
the truncation as the following examples show.

\begin{exa}\label{4.3} The area of the triangles in Fig. 5.4
and 5.5 are $a\cdot i=\pi-(\frac{\pi}2+\frac{\pi}2-ia)$ and
$-A=\pi-(\frac{\pi}2+\frac{\pi}2+A)$ respectively.
\end{exa}

\begin{figure}[h]
\begin{center}
\includegraphics[width=0.28\textwidth]{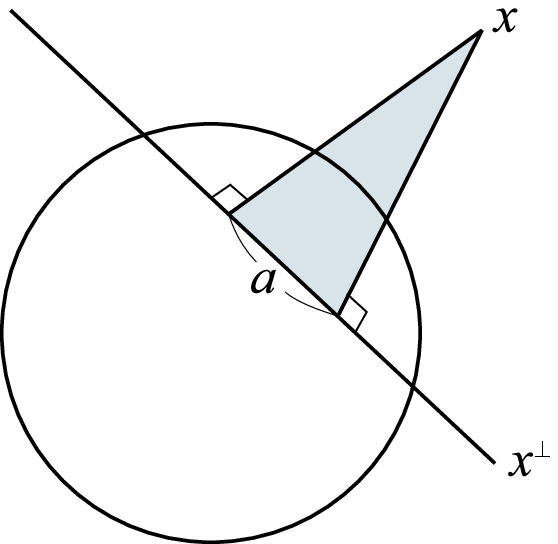}\quad\quad\quad
\includegraphics[width=0.35\textwidth]{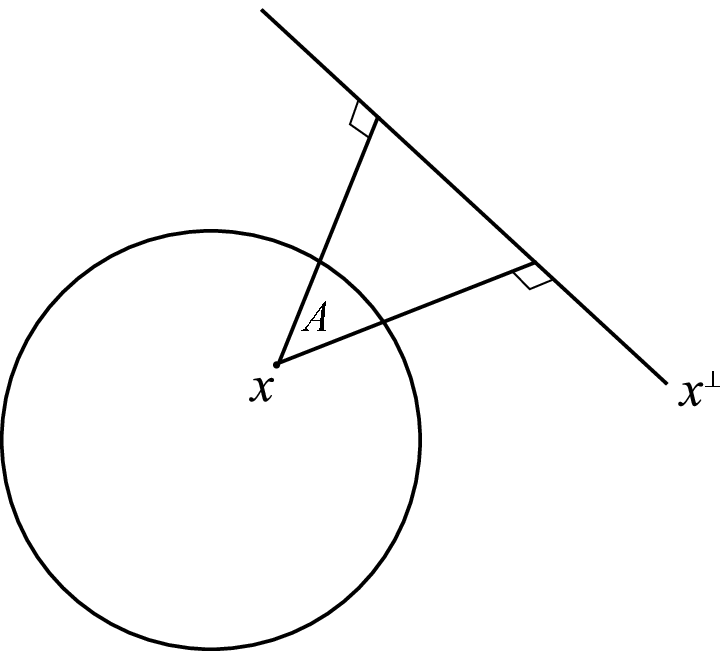}
\caption{\textbf{5.4}\textbf{\hskip 5truecm Fig 5.5\quad}}
\end{center}
\end{figure}

\begin{exa}\label{4.4} The triangle (1,2,3) in Fig 5.6 has
area $\pi-(A-bi-ei)= \pi-A+(b+e)i$. Also the triangle (1,2,3) is
divided into three polygons (2,4,5), (3,6,7) and (1,4,5,7,6) and
each polygon has pure imaginary area $ei,bi$ and real value area
$\pi-A$ respectively.
\end{exa}

\begin{figure}[h]
\begin{center}
\includegraphics[width=0.22\textwidth]{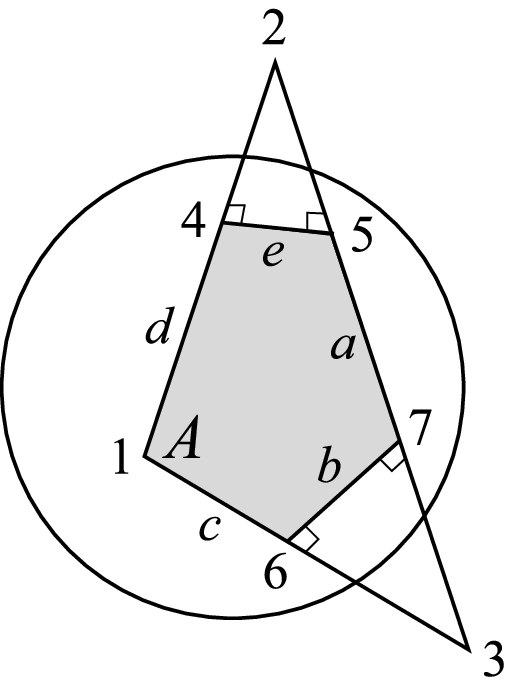}\quad\quad\quad\quad
\includegraphics[width=0.32\textwidth]{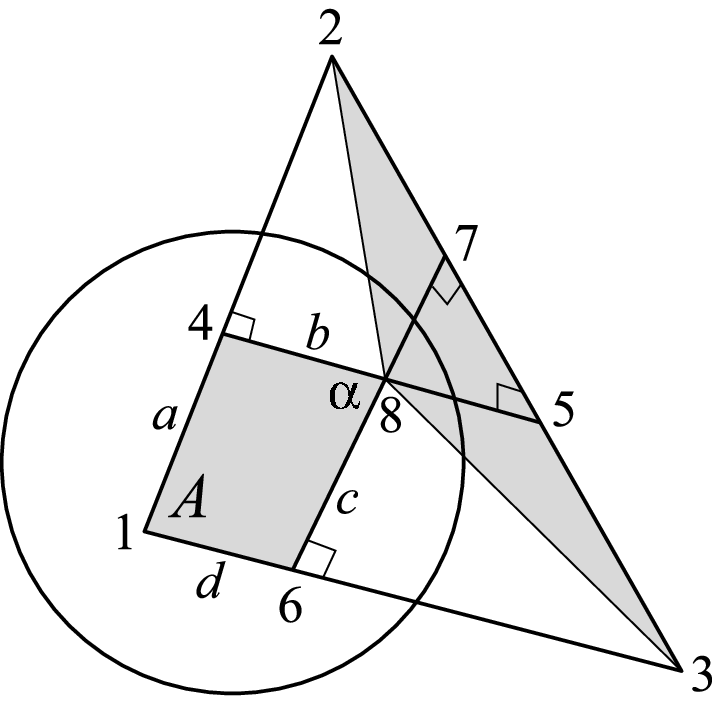}
\caption{\textbf{5.6}\textbf{\hskip 5truecm  Fig 5.7\quad\quad}}
\end{center}
\end{figure}

\begin{exa}\label{4.5} The triangle (1,2,3) in Fig 5.7 is
divided into various polygons. Notice that the line segment (6,7)
is $3^{\bot}$ and (4,5) is $2^{\bot}$, and hence (2,3) is
$8^{\bot}$ which implies that the length of the line segment (5,8)
is equal to $\frac{\pi}2 i$ and similarly for (7,8). Let $b$ be
the length of (4,8), $c$ for (6,8) and $\alpha$ be the angle
$\angle(4,8,6)$. Then we have
$$
\aligned
    \Area (1,2,3) &= \pi-(A+({\pi \over 2}-bi)
                    +({\pi\over 2} - ci))\\
                &=-A+(b+c)i,
\endaligned
$$
and
$$
\aligned
    \Area (1,4,8,6) + &\Area (2,8,3) +
        \Area (2,4,8) + \Area (3,6,8) \\
    &=(\pi-A-\alpha) +(\alpha - \pi) +(bi) +(ci)\\
    &=-A+(b+c)i
\endaligned
$$
Notice that the real part is the area of the shaded quadrangle and
the imaginary part is the total length of the truncated sides of the
quadrangle.
\end{exa}

\begin{exa}\label{4.6} From the outer triangle which is dual to the inner triangle $\triangle(1,2,3)$
 in Fig 5.8, we get
$$
    \Area (1',2',3') = \pi - (\pi-ai+\pi - bi +\pi - ci)
            = -2\pi +(a+b+c)i.
$$
If the inner triangle in the figure shrinks to one point, then the
area of the triangle becomes $-2\pi (=\text{vol }(\Bbb S^2_H)/2)$.
In fact in this case, the vertices of the triangle lies on the
equator and the angle at the vertices become $\pi$.
\end{exa}

\begin{figure}[h]
\begin{center}
\includegraphics[width=0.45\textwidth]{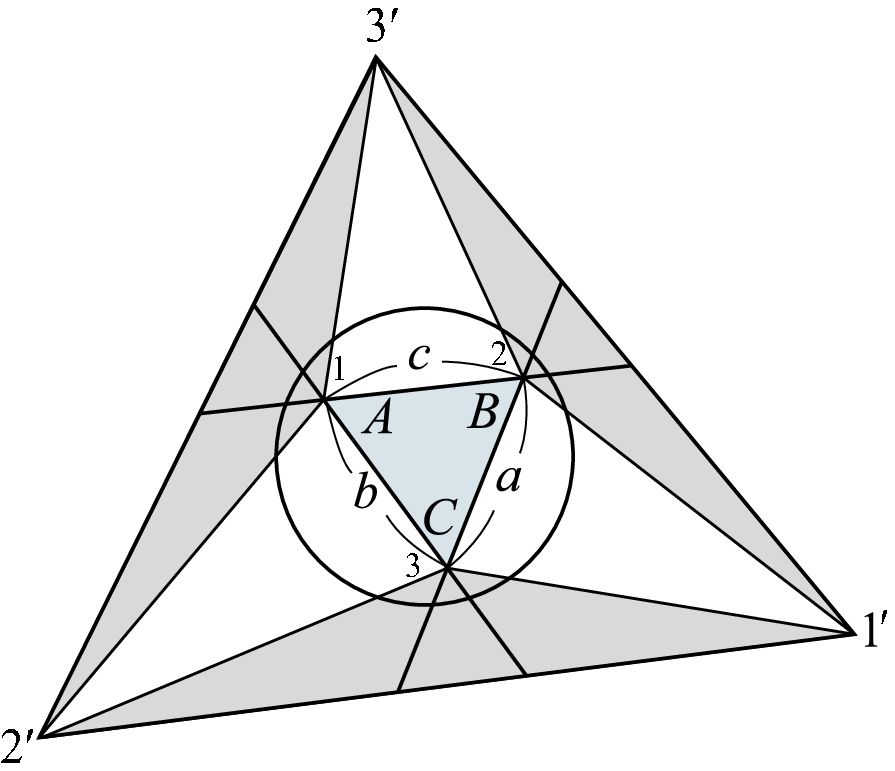}
\caption{\textbf{5.8}}
\end{center}
\end{figure}

We can generalize hyperbolic trigonometry to any triangle in the
extended hyperbolic space $\Bbb S^2_H$ or in general $\Bbb S^n_H$
(see \cite{C}). Naturally, already known cosine and sine laws at the
hyperbolic space can be applicable to the extended hyperbolic space.
Then all the known trigonometric formulas for special hyperbolic
polygons are in fact explained in a consistent manner through the
trigonometry for various triangles of $\Bbb S^2_H$ lying across the
ideal boundary. See Fig. 5.6 and 5.7 for the right angled polygons
and the corresponding triangles in $\Bbb S^2_H$. All the cosine and
sine laws of quadrangles in Fig. 5.6 and 5.7 could be obtained from
the cosine and sine laws of the outer triangles $\triangle(1,2,3)$
in the figures. For example, the following cosine law of the
pentagon in Fig. 5.6 $\cosh a=\frac{\cosh b\cosh e+\cos A}{\sinh b
\sinh e}$ is deduced from the following generalized cosine law for
the triangle $\triangle(1,2,3)$ in $\Bbb S^2_H$,
$$\cosh (a+\pi i)=\frac{\cos (-bi)\cos (-ei)+\cos A}{\sin (-bi) \sin
(-ei)}.$$ 

Also the extended hyperbolic space shows an evident and
geometric reason of resemblance of laws between the hyperbolic space
and spherical space. Furthermore we can as well derive the trigonometry for the
Lorentzian triangles and polygons, which has not been determined
well as far as we know. We refer the interested reader to see \cite{C} for the details.

\vskip 1pc
From Theorem \ref{4.1}, we can also induce the volume of elementary
three dimensional objects, such as lens and trihedron $\Tri (A,B,C)$
with dihedral angles $A,B,C$.

\begin{cor}\label{4.7} The $3$-dimensional volume of a convex lens
$l(A)$ with dihedral angle $A$ is $-A\pi i$.
\end{cor}
\begin{proof}
$$
\aligned &{\vol (\Bbb S^3_H) \over \vol (\Bbb S^2_H)} \times
    \text{Area of 2-dimensional cross section for the lens of angle }A \\
=& {\pi \over 2}i\times (-2A).
\endaligned
$$
\end{proof}

\begin{cor}\label{4.8} The volume of a trihedron $\Tri (A,B,C)$
with dihedral angles $A,B,C$ is represented as
$$
\vol (\Tri(A,B,C)) = {\pi\over 2} i (\pi -A-B-C).
$$
\end{cor}

The volume of a hyperbolic trihedron has a similar expression as the
spherical one which has a well-known volume formula, $V={\pi \over
2}(A+B+C-\pi)$. The difference between the two formulas is $i^3$.

\begin{exa}\label{5.11}
Consider a tetrahedron $\triangle^3$ some of whose vertices (but
not edges) are lying outside the hyperbolic part. Let $P$ be a
hyperbolic polyhedron truncated from $\triangle^3$ by the dual
planes of outside vertices and $T$ be the truncated part of
$\triangle^3$ so that $\triangle^3=P\cup T$. See Fig. 5.9. Then we
have
$$\al
\text{vol }(\triangle^3)&=Re~\text{vol }(\triangle^3)+Im~\text{vol
}(\triangle^3)\\&=\text{vol }(P)+\text{vol }(T)\\&=\text{vol
}(P)+\frac{\pi}{4}i\times (\text{area of the truncated face}).\eal
$$
This follows immediately from Theorem 5.1 since $T$ is exactly the
half of the lune containing $T$.

\begin{figure}[h]
\begin{center}
\includegraphics[width=0.25\textwidth]{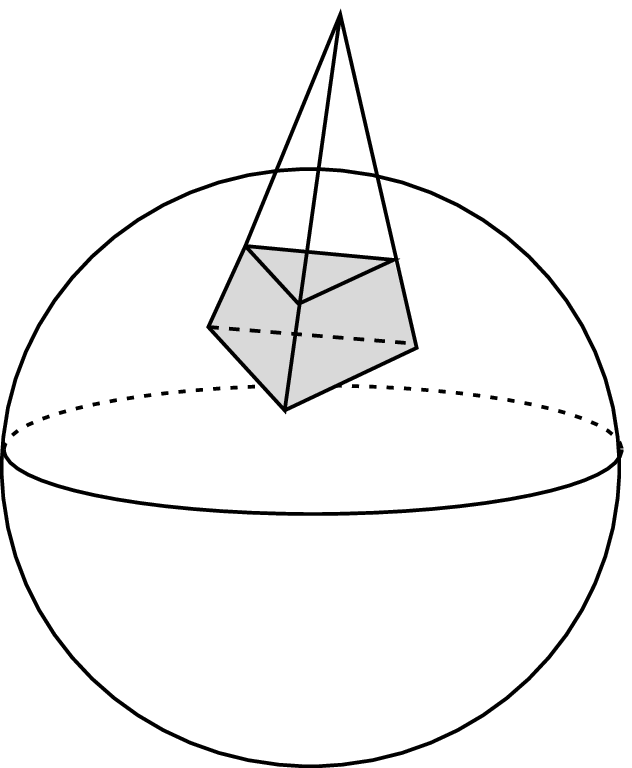}
\caption{\textbf{5.9}}
\end{center}
\end{figure}

Ushijima \cite{U} found that $\text{vol }(P)$ is the real part of
the value which is obtained by applying Murakami-Yano's volume
formula for hyperbolic tetrahedron \cite{MY} formally to
$\triangle^3$. In this case we have a problem of determining a
"right" value among the multi-values from the formula and in our
extended model the above $\text{vol }(\triangle^3)$ is the correct
value and the pure imaginary part is interpreted as the volume of
$T$ as well as the area of the truncated face.

The general statement for such phenomenon requires a proof of the
analyticity of the volume formula for tetrahedra across the
boundary in the extended model and will be deferred to a
subsequent paper.
\end{exa}

We consider another interesting application of the extended model.
Let $\triangle (1,2,3)$ be a convex spherical triangle with
vertices 1, 2, 3 on the 2-sphere $\Bbb S^2\subset \Bbb R^3$ and
$\triangle (1',2',3')$ be its dual triangle so that the distance
between a vertex $i'$ ($i$, resp.) and each point on the edge
$(j,k)$ ($(j',k')$, resp.) is $\pi\over 2$, i.e., they are
perpendicular as vectors in $\Bbb R^3$, where
$\{i,j,k\}=\{1,2,3\}$. Then $\Bbb S^2$ can be decomposed into 8
triangles, $\triangle (1,2,3)$, $\triangle (1',2,3)$, $\triangle
(1,2',3)$, $\triangle (1,2,3')$, $\triangle (1',2',3)$, $\triangle
(1',2,3')$, $\triangle (1,2',3')$, $\triangle (1',2',3')$. If we
denote the volume of $\triangle (i,j,k)$ again by $\triangle
(i,j,k)$ abusing the notation, then we have
$$
\al &\triangle (1,2,3)+\triangle (1',2',3)+\triangle
(1',2,3')+\triangle (1,2',3')\\
=&(A+B+C-\pi)+(\pi-C+{\pi\over 2}+{\pi\over 2}-\pi)+(\pi-B+{\pi\over 2}+{\pi\over 2}-\pi)+
(\pi-A+{\pi\over 2}+{\pi\over 2}-\pi)\\=&2\pi,\\
\eal$$

$$\al
&\triangle (1',2,3)+\triangle (1,2',3)+\triangle
(1,2,3')+\triangle
(1',2',3')\\
=&(a+{\pi\over 2}+{\pi\over 2}-\pi)+(b+{\pi\over 2}+{\pi\over
2}-\pi)+(c+{\pi\over 2}+{\pi\over
2}-\pi)+(\pi-a+\pi-b+\pi-c-\pi)\\=&2\pi, \eal
$$
where $A,B,C$ are the angles at vertices 1,2,3 respectively and
$a,b,c$ are there opposite edge lengths, which are the same as the
angles at dual vertices $1',2',3'$ respectively.

We can express these identities equivalently as a single one,
denoting $n'$ by $n^{-1}$, as follows.
$$\sum_{\epsilon_1,\epsilon_2,\epsilon_3=\pm
1}\epsilon_1\epsilon_2\epsilon_3~\triangle(1^{\epsilon_1},2^{\epsilon_2},3^{\epsilon_3})=0$$

We can generalize these identities to a convex polygon as follows:
$$\text{vol }(\rm I)+\text{vol }(\rm III)=\text{vol }(\rm II)+\text{vol
}(\rm IV)=2\pi$$ or equivalently,
$$\text{vol }(\rm I)-\text{vol }(\rm II)+\text{vol }(\rm III)-\text{vol
}(\rm IV)=0,$$ where the regions $I$ to $IV$ are depicted in Fig
5.10.

\begin{figure}[h]
\begin{center}
\includegraphics[width=0.45\textwidth]{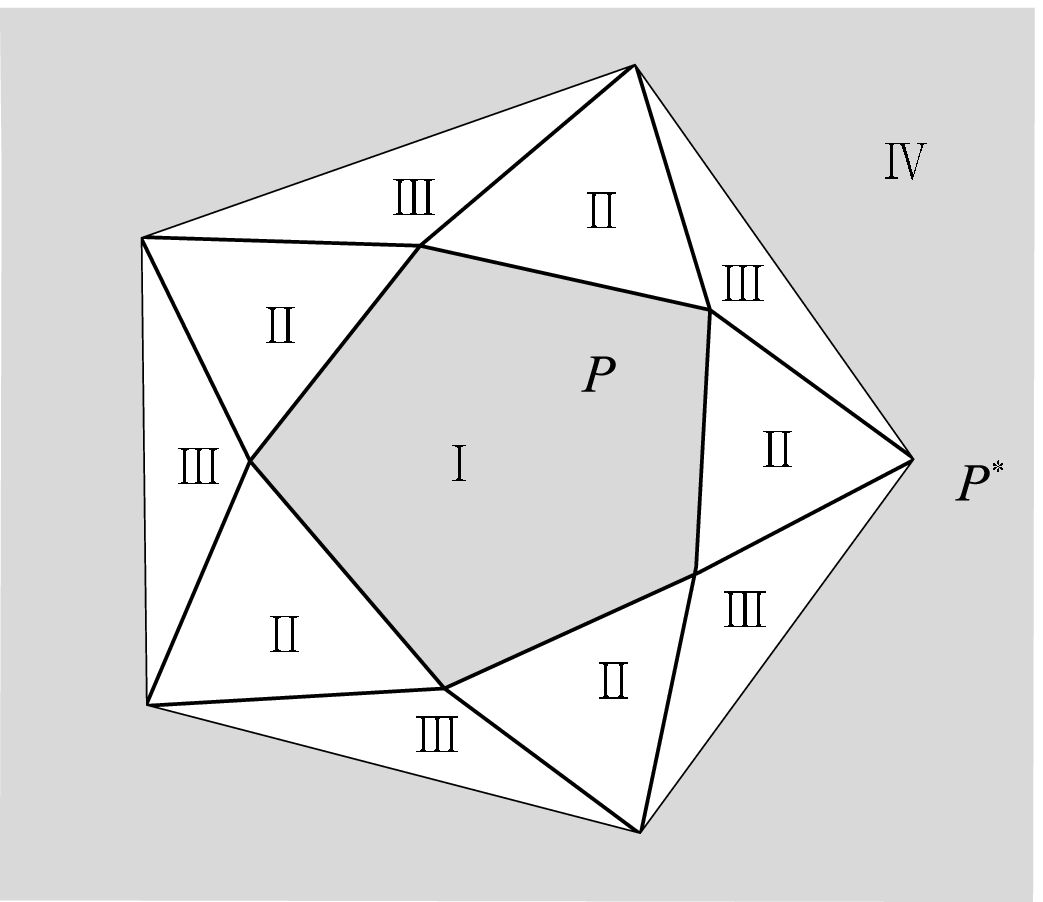}
\caption{\textbf{5.10}}
\end{center}
\end{figure}

Let's generalize these observations to the $n$-dimensional case.
First start with a convex polyhedron $P$ and its dual $P^*$ on
$\Bbb S^n$. Let $P$ be the region $\rm I$ and consider the the
region $\rm II$ consisting of polyhedra having ($n-1$)-dimensional
face in common with ${\rm I} =P$, and the region $\rm III$ of
polyhedra having ($n-1$)-dimensional face in common with $\rm II$,
$\cdots$, and so on, so that we obtain $P^*$ as the ($n+2$)-th
region. Then we can expect the following identity as above.
$$\text{vol }({\rm I})-\text{vol }({\rm II})+\text{vol }({\rm III})-\cdots+(-1)^{n+1}\text{vol
}(P^*)=0.$$ When $n=3$, this becomes
$$\text{vol }(\rm I)-\text{vol }(\rm II)+\text{vol }(\rm III)-\text{vol
}(\rm IV)+\text{vol }(V)=0,$$ or equivalently,

\begin{align}
                   &\text{vol }(P)+\text{vol }({\rm III})+\text{vol }(P^*)=\pi^2\label{sm}\\
                   &\text{vol }({\rm II}) +\text{vol }({\rm IV})=\pi^2.\nonumber
\end{align}
Indeed the identity (\ref{sm}) is known as Santal\'o-Milnor
relation, and is proved by Santal\'o using integral geometry for a
3-simplex \cite{S} and proved by Milnor using Schl\"afli formula
for a general convex polyhedron \cite{M}. This also can be proved
using elementary geometry \cite{C2}.

We can obtain a simplex result for a hyperbolic convex polyhedron
on $\Bbb H^n$ using the hyperbolic sphere $\Bbb S^n_H$. Consider a
hyperbolic triangle $\triangle (1,2,3)$ and its dual $\triangle
(1',2',3')$ on $\Bbb S^2_H$. Then as for the spherical case, we
obtain
\begin{equation}\label{ep}\sum_{\epsilon_1,\epsilon_2,\epsilon_3=\pm
1}\epsilon_1\epsilon_2\epsilon_3~\triangle(1^{\epsilon_1},2^{\epsilon_2},3^{\epsilon_3})=0,
\end{equation} where $1^{-1}=1',\ldots$ etc.
Here the dual triangle $\triangle (1',2',3')$ is not a convex hull
of its vertices, but rather its complement to obtain a
decomposition of $\Bbb S^2_H$ into 8 triangles as before. Also we
consider only a {\it compact} hyperbolic triangle $\triangle
(1,2,3)$ in the formula (\ref{ep}). In fact, (\ref{ep}) does not
hold for a more general type triangle $\triangle (1,2,3)$, even
for an ideal triangle. In these cases we have to correct
(\ref{ep}) by some constants.

For a compact hyperbolic convex polygon, the following identities
are obtained easily as for the spherical case.
$$\text{vol }({\rm I})+\text{vol }({\rm III})=\text{vol }({\rm II})+\text{vol
}({\rm IV})=\frac{\text{vol }(\Bbb S^2_H)}{2}$$ or equivalently,
$$\text{vol }(\rm I)-\text{vol }(\rm II)+\text{vol }(\rm III)-\text{vol
}(\rm IV)=0.$$ See Fig 5.10 for the description of regions.

Again for $n$-dimensional case, we can expect the following
identity for a hyperbolic convex polyhedron $P$ and its dual
$P^*$, as the complement of the convex hull of its vertices in
$\Bbb S^n_H$.

\begin{equation}\label{S}\text{vol }(P)-\text{vol }({\rm II})+\text{vol }({\rm
III})-\cdots+(-1)^{n+1}\text{vol }(P^*)=0.\end{equation} When
$n=3$, we have
$$\text{vol }(\rm I)-\text{vol }(\rm II)+\text{vol }(\rm III)-\text{vol
}(\rm IV)+\text{vol }(V)=0$$ or equivalently,

\begin{align}
                   &\text{vol }(P)+\text{vol }({\rm III})+\text{vol }(P^*)=\frac{\text{vol }(\Bbb S^3_H)}{2}
                   \label{sm'}\\
                   &\text{vol }({\rm II}) +\text{vol }({\rm IV})=\frac{\text{vol }(\Bbb
                   S^3_H)}{2}.\nonumber
\end{align}

The identity (\ref{sm'}) can be rewritten as
\begin{equation}\label{sm''}
                   \text{vol }(P)+\text{vol }({\rm III})+\text{vol
                   }(P')=0,
\end{equation}
where $P'=P^*\cap K^3$.

If $P$ has an ideal vertex, $\text{vol }({\rm III})$ and
$\text{vol }(P')$ become $-\infty$ and $\infty$ respectively. The
formula (\ref{sm''}) is proved by Santal\'o using integral
geometry for a compact hyperbolic 3-simplex \cite{S}, proved by
Su\'arez-Peir\'o using Schl\"afli formula for a compact hyperbolic
simplex in general dimension \cite{E}, and by J. Murakami and
Ushijima using a volume formula for a hyperbolic simplex \cite{J}.
Again this can also be proved by elementary geometry on $\Bbb
S^3_H$ using Theorem 5.1 \cite{C2}.

The term $\text{vol }({\rm III})$ in (\ref{sm''}) is the sum of
the products of each edge length and its dihedral angle and
interpreted as mean curvature, but in the setting of extended
model $\Bbb S^3_H$, it simply represents the volume of tetrahedra
in  between $P$ and $P'$ lying across the ideal boundary
$\partial\Bbb H^3$. Similarly the terms $\text{vol }({\rm II}),
\text{vol }({\rm III}), \ldots,$ so called integral mean
curvatures, are interpreted as the volumes of collections of
tetrahedra sitting in between $P$ and $P'$ successively in the
extended model $\Bbb S^n_H$.

\centerline{\bf APPENDIX} \vskip 1pc {\it Let $L$ be a "wedge
region" contained in the Lorentzian part $\Bbb L^n$ determined by
two hyperplanes $\mathcal P_1$ and $\mathcal P_2$ such that
$\mathcal P_1\cap \mathcal P_2$ is tangent to $\partial\Bbb H^n$
(see the picture below). The volume of $L$ is finite if none of
$\mathcal P_i, i=1,2$. is tangent to $\partial\Bbb H^n$, and is
infinite otherwise. When $n=2$, we already know this. Or it
follows from direct calculation similar the following integrals.}

\begin{figure}[h]
\begin{center}
\includegraphics[width=0.6\textwidth]{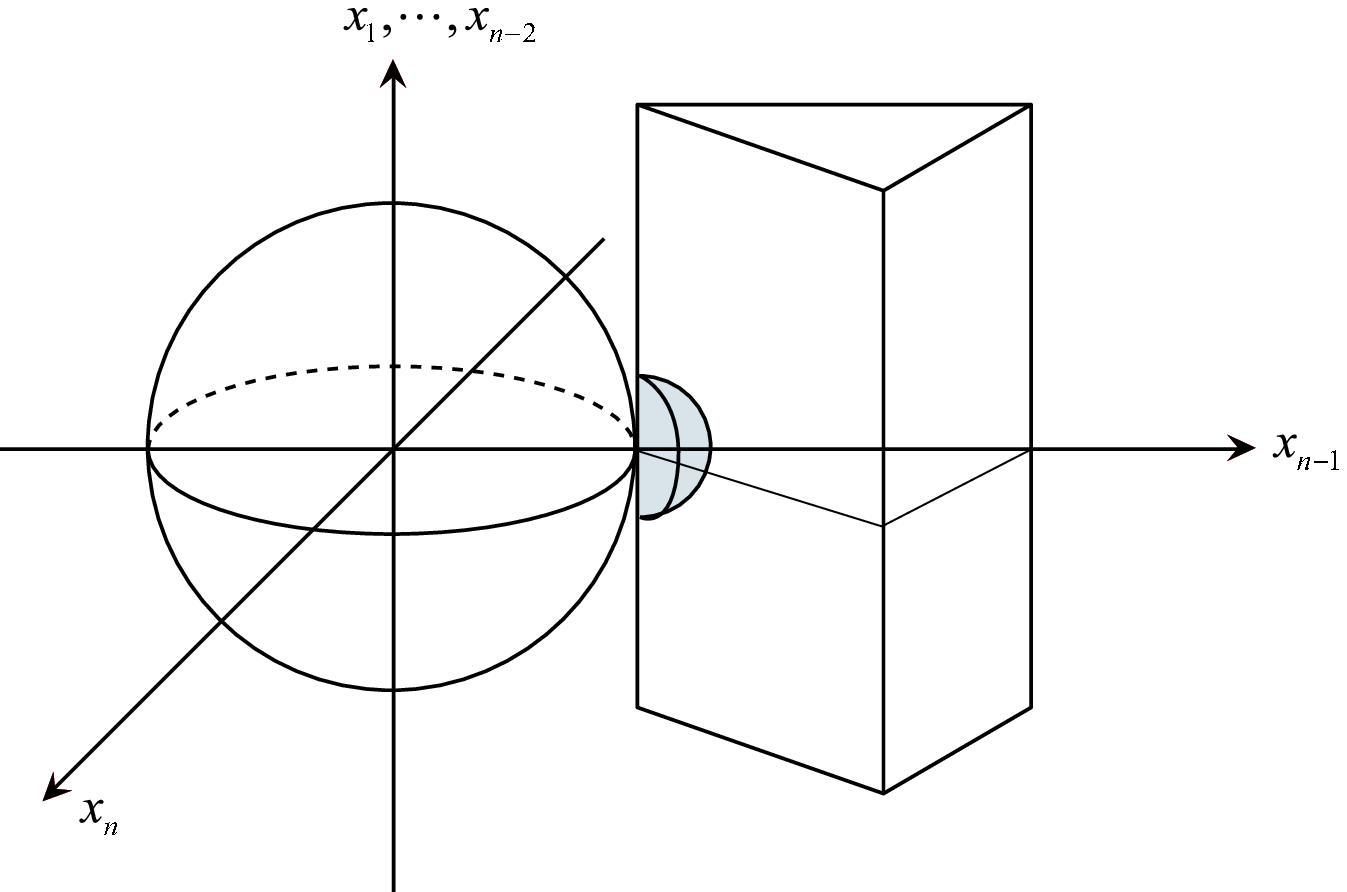}
\end{center}
\end{figure}

\begin{proof}
When $n\ge 3$, it suffices to show that the integral
$$I=\int^{\infty}_{-\infty}\cdots\int^{\infty}_{-\infty}\int^{\infty}_{1}\int^{k(x_{n-1}-1)}_0
\frac{ dx_n dx_{n-1} dx_{n-2}\cdots
dx_1}{(x_1^2+\cdots+x_n^2-1)^{\frac{n+1}{2}}}$$ is finite for all
$k>0$ and is infinite for $k=\infty$.

We will denote $\int f\sim \int g$ when $\int f<\infty$ iff $\int
g<\infty$. Consider a spherical coordinate
$(\rho,\theta_1,\ldots,\theta_{n-1})$ centered at the tangent
point $(0,\ldots,0,1,0)$:
$$\al
x_1&=\rho\cos \theta_1\\
x_2&=\rho\sin \theta_1\cos \theta_2\\
\vdots\\
x_{n-1}&=\rho\sin \theta_1\sin \theta_2\cdots\sin\theta_{n-2}\cos \theta_{n-1}\\
x_{n}&=\rho\sin \theta_1\sin \theta_2\cdots\sin\theta_{n-2}\sin \theta_{n-1},\\
\eal$$
 where $\rho=|x|$, $\theta_i=\angle(e_i, x_i e_i+\cdots+x_n
e_n)$ with $0\le \theta_i\le\pi$ if $i<n-1$, and $\theta_{n-1}$ is
the polar angle from $e_{n-1}$ to $x_{n-1} e_{n-1}+x_n e_n$ with
$0\le \theta_i\le 2\pi$ (see \cite[p.45]{Ra}). Using the spherical
coordinate with $\alpha=\tan^{-1} k$, $0<\alpha<\frac{\pi}{2}$,
the integral $I$ becomes as follows.
$$\al
I&=\int^{\alpha}_0\int^{\pi}_0\cdots\int^{\pi}_0\int^{\infty}_0\frac{\rho^{\frac{n-3}{2}}\sin^{n-2}\theta_1
\sin^{n-3}\theta_2\cdots\sin\theta_{n-2}}{(\rho+2\sin\theta_1
\sin\theta_2\cdots\sin\theta_{n-2}\cos\theta_{n-1}
)^{\frac{n+1}{2}}}~~d\rho
d\theta_1\cdots d\theta_{n-1}\\
&\sim\int^{\alpha}_0\int^{\pi}_0\cdots\int^{\pi}_0\int^{\delta}_0\frac{\rho^{\frac{n-3}{2}}\sin^{n-2}\theta_1
\sin^{n-3}\theta_2\cdots\sin\theta_{n-2}}{(\rho+2\sin\theta_1
\sin\theta_2\cdots\sin\theta_{n-2}\cos\theta_{n-1}
)^{\frac{n+1}{2}}}~~d\rho d\theta_1\cdots d\theta_{n-1}=:I', \eal
$$
where $\delta>0$ and notice that the region $\rho\ge\delta$ being
a relatively compact region inside $\Bbb L^n$ has obviously finite
volume. Integrate with respect to $\rho$ using the formula
$$\int \frac{x^n}{(x+a)^{n+2}} ~~dx=\frac{1}{(n+1)a}\frac{x^{n+1}}{(x+a)^{n+1}},
$$
we have
$$\al
I'&=\int^{\alpha}_0\int^{\pi}_0\cdots\int^{\pi}_0\frac{\sin^{n-2}\theta_1
\sin^{n-3}\theta_2\cdots\sin\theta_{n-2}}{(n-1)\sin\theta_1
\cdots\sin\theta_{n-2}\cos\theta_{n-1}}\frac{\delta^{\frac{n-1}{2}}~d\theta_1\cdots
d\theta_{n-1}}{(\delta+2\sin\theta_1
\cdots\sin\theta_{n-2}\cos\theta_{n-1} )^{\frac{n-1}{2}}}\\
&\sim\int^{\alpha}_0\int^{\pi}_0\cdots\int^{\pi}_0\frac{\sin^{n-3}\theta_1
\sin^{n-4}\theta_2\cdots\sin\theta_{n-3}}{\cos\theta_{n-1}}~~d\theta_1\cdots
d\theta_{n-1}\\
&\sim\int^{\alpha}_0 \sec\theta_{n-1}~~d\theta_{n-1}.\eal$$
 Now this final integral is clearly finite if
$\alpha<\frac{\pi}{2}$, and is infinite if $\alpha=\frac{\pi}{2}$,
i.e., $k=\infty$.
\end{proof}

{\sc \begin{flushleft} Department of Mathematics, University of Seoul, Seoul 130-743,
Korea\\[5pt]
  Department of Mathematics, Seoul National University, Seoul 151-742, Korea\\[5pt]
E-mail: yhcho@uos.ac.kr\\
\mbox{\phantom{E-mail: }hyukkim@snu.ac.kr}
\end{flushleft}}
\end{document}